\documentclass[pdflatex,sn-mathphys-num]{sn-jnl}

\usepackage{graphicx}
\let\orcidlogo\relax

\usepackage{orcidlink}
\usepackage{amsmath,amssymb}
\usepackage{bm}
\usepackage{booktabs}
\usepackage{subcaption}
\usepackage{adjustbox}
\usepackage{algorithm}
\usepackage{algpseudocode}

\theoremstyle{definition}

\newcommand{\Bp}[2]{\mathcal{B}^{#1,+}_{#2}}
\newcommand{\Bm}[2]{\mathcal{B}^{#1,-}_{#2}}
\newcommand{\dd}{\mathrm{d}}

\begin{document}

\title[An Active Flux Method for the three-dimensional Maxwell's equations]{An Active Flux Method for Maxwell's equations}

\author[1]{\fnm{Yanick} \sur{Kiechle}\orcidlink{0009-0002-8529-7482}}\email{yanick.kiechle@hhu.de}
\author[1]{\fnm{Christiane} \sur{Helzel}\orcidlink{0000-0002-9901-7149}}\email{christiane.helzel@hhu.de}

\affil*[1]{\orgdiv{Faculty of Mathematics and Natural Sciences, Institute of Mathematics}, \orgname{Heinrich-Heine-University D\"usseldorf}, \orgaddress{\city{D\"usseldorf}, \country{Germany}}}

\abstract{We present a novel Active Flux method and a divergence-free method using exact evolution operators for solving Maxwell's equations in three space dimensions on Cartesian grids.
	The proposed methods are fully discrete, third-order accurate and rely on a truly multidimensional evolution operator 
	based on Kirchhoff's formula for the electromagnetic wave equation. 
	The first method is a classical Active Flux method with globally continuous reconstruction while the second 
	method is an adaptation of the Active Flux method to a locally divergence-free reconstruction by Balsara.
	The applied reconstruction is built from facial degrees of freedom and is no longer globally continuous. 
	The classical method keeps the divergence error bounded.
	The divergence-free method with face averages initialized by a potential keeps the constraints at round-off. We compare the accuracy of both methods 
	for smooth and discontinuous test problems. }


\maketitle
{\renewcommand{\thefootnote}{}\footnotetext{Funded by the Deutsche Forschungsgemeinschaft (DFG, German Research Foundation) through FOR 5409 \textit{Structure-preserving Numerical Methods for Bulk and Interface Coupling of Heterogeneous Models}, within the project \textit{An Active Flux Method for the Vlasov--Maxwell System}, project number 463312734.}}

\section{Introduction}

In this publication, we present the first fully discrete Active Flux method for Maxwell's equations in three dimensions.
In 2011, Roe and his former student Eymann introduced Active Flux methods for advection and acoustics in one and two dimensions \cite{AFSchemesEymannRoe, MultiDAFSchemesEymannRoe}.
For the simplest case of one-dimensional advection, the method agrees with van Leer's method V from \cite{VANLEER1977276}.
Its defining feature is a continuous reconstruction, which removes the need for Riemann solvers and allows the computation of numerical fluxes via Simpson's rule.
Roe later argued that third order represents a suitable candidate, being an odd-order method \cite{tyacke_musings_2025} with not too much complexity \cite{roe_is_2017}.
In one and two spatial dimensions, higher-order approximations have been considered \cite{abgrall_extensions_2023,barsukow_generalized_2026}. In this publication
we restrict ourselves to a third-order accurate method.
Besides the cell-average values, Active Flux uses point values as additional degrees of freedom. 
These define the reconstruction and enter the numerical fluxes.
The point values allow for a compact stencil in space and time while retaining third-order accuracy. 
This makes Active Flux methods attractive to parallelization and capable of achieving high accuracy on coarse grids owing to their low numerical dissipation.
Interest in this method grew, with extensions to Cartesian grids \cite{helzel_new_2019, barsukow_active_2019,CHUDZIK2021125501}. 
With appropriate evolution operators for the point values,
Active Flux presents a finite volume approach that tries to incorporate the true physics
and achieves a numerical domain of dependence in close agreement with the physical one.
Recent publications continue to address this matter by developing evolution operators for
the Euler equations \cite{chudzik_active_2024,chudzik_fully_2025,barsukow_active_2025,helzel2026fullydiscreteactiveflux}. 
Although Active Flux was initially formulated as a fully discrete method, researchers have since developed semi-discrete methods. 
Among these are the method of Barsukow, Klingenberg and coworkers and their approach to the 
\textit{generalized} Active Flux methods \cite{barsukow_generalized_2026, Duan_Barsukow_Klingenberg_2025, abgrall2025semi}, which have also been applied to the ideal MHD equations \cite{duan2025activefluxidealmagnetohydrodynamics}, as well as 
Abgrall et al. who developed the PAMPA schemes that use the same underlying degrees of freedom \cite{abgrall_extensions_2023,abgrall_novel_2025, ABGRALL2026114521}.
Recently, works have been published that presented Active Flux methods for the Vlasov--Poisson 
and Vlasov--Maxwell systems \cite{HENSEL2025114294,KIECHLE2025113693, gruenwald2025}. 
In kinetic plasma models such as the Vlasov--Maxwell system, Maxwell's equations describe the propagation of electromagnetic waves coupled to the Vlasov equation, which describes the evolution of a plasma~\cite{Kraus_Kormann_Morrison_Sonnendrücker_2017}.
In \cite{gruenwald2025} the newly developed Active Flux methods for the Vlasov equation were combined with
a classical Yee scheme to provide a new method for the Vlasov--Maxwell system. 
The authors showed that Active Flux methods for the Vlasov equation have significantly lower dissipation than the conservative PFC method \cite{FILBET2001166}.
This motivates our work, where we propose an Active Flux method for Maxwell's equations in three dimensions.
The main difficulty with such methods is constructing a suitable evolution operator. This can be resolved exactly using a multidimensional operator based on Kirchhoff's formula. 
Combined with Active Flux's low numerical dissipation, this promises high accuracy on coarse grids without the need for Riemann solvers or a separate time integration method.

The development of numerical methods for Maxwell's equations has a long history. The most prominent method is the Yee scheme \cite{kane_yee_numerical_1966, TAFLOVE2005629}, a second-order 
finite difference time-domain method employing a staggered grid in space and time. Many subsequent methods mimic the structure of the de~Rham complex, starting from edge elements \cite{nedelec_mixed_1980}. 
They evolve differential forms in function spaces chosen to respect the exterior calculus, following Arnold, Falk and Winther 
on finite element exterior calculus~\cite{Arnold_Falk_Winther_2006}. 
These structure-preserving discretizations include both finite difference \cite{bach2025sbpfdecsummationbypartsfinitedifference} and finite element approaches \cite{KORMANN2025118290}, but are typically 
semi-discrete in time, requiring a suitable time integration method.
Alternatively, the divergence constraints can be controlled by divergence cleaning \cite{MUNZ2000484} or enforced locally by divergence-free basis functions in discontinuous Galerkin methods \cite{COCKBURN2004588}.

In the present work we propose two fully discrete methods for Maxwell's equations
without source terms in three spatial dimensions. The first uses a tensor-product construction from
lower dimensions, resulting in $27$ degrees of freedom per cell for the second-degree
polynomial reconstruction required for third-order accuracy. The time step is restricted by
a CFL condition of $0.5$ to maintain the compact stencil, and point values are evolved
using a truly multidimensional exact evolution operator for the homogeneous Maxwell system.
The evolution operator is derived from the wave equation and is based on Kirchhoff's formula. 
This method preserves the divergence constraints only up to the truncation error of its reconstruction.
The resulting divergence error stays bounded but does not vanish. This method is an extension 
of Active Flux methods to Maxwell's equations in three dimensions. It relies on the previous ideas 
of evolution operators for wave dominated systems \cite{eymann_multidimensional_2013, lukacova-medvidova_evolution_2004, chudzik_active_2024, chudzik_fully_2025, barsukow_active_2025} and 
extends the ideas of Active Flux methods for Maxwell's equations in two dimensions \cite{chudzik_active_2024} to three dimensions.
In 2D, Maxwell's equations reduce to two sets of acoustic equations, where a discrete divergence is preserved by an exact evolution operator.

As the second method, we propose one that relies on the same evolution operator but has a reconstruction 
that is locally divergence-free. The reconstruction is based on the one derived by Balsara in \cite{BALSARA2001614,BALSARA2004,BALSARA2009} for MHD 
and applied to Maxwell's equations in \cite{HAZRA2019298,BALSARA20191108}.
In that setting, the reconstruction is designed to be divergence-free for the magnetic field. In our proposed 
method, we adapt this reconstruction to both the electric and magnetic field, to obtain a divergence-free method 
for Maxwell's equations. In \cite{BALSARA2009} the reconstruction is 
built from facial moments, which are reconstructed from the face averages of neighboring faces.
In our approach, we compute these remaining degrees of freedom from point values on the cell faces, which are 
evolved using the same evolution operator as in the classical Active Flux method. 
Since the reconstruction is not continuous, the curl in the evolution operator is understood in the distributional sense.
We use appropriate finite difference methods to compute the facial moments from the face point values.

We update the face averages of $\bm B$ and $\bm E$ via Stokes' theorem to maintain the divergence at machine precision. 
The arising edge integrals are approximated using appropriate quadrature rules from the evolved point values. 
The reconstruction is not continuous and does not interpolate point values on the boundaries.
Consequently, point values on the boundary of a cell are not uniquely defined at time $t_n$. Therefore, we replace Simpson's rule 
as the quadrature rule in time by a right Radau rule. 
This gives up the continuity of the reconstruction for the divergence-free property. 

The article is structured as follows: Section~\ref{sec:AF} introduces the Active Flux
framework. Section~\ref{sec:Maxwell} derives the evolution operators for
Maxwell's equations without source terms. Section~\ref{sec:divfree} presents the divergence-free
method. Numerical experiments are presented in
Section~\ref{sec:numerics}, and Section~\ref{sec:conclusion} concludes this article.

\section{Fully discrete Active Flux methods in three dimensions}\label{sec:AF}
To illustrate the formulation of Active Flux, we will derive the general form of the fully discrete method for a general 
three-dimensional conservation law on a domain $\Omega \subset \mathbb{R}^3$. Consider the vector of conserved quantities $\bm q:\Omega \times \mathbb{R}^+ \to \mathbb{R}^m$ with a flux function $\bm{f}$
defined by $\bm{f}(\bm q) := (f(\bm q), g(\bm q), h(\bm q)) : \mathbb{R}^m \to \mathbb{R}^{3 \times m}$. Then the general
conservation law is given by the following partial differential equation
\begin{align}
		\frac{\partial}{ \partial t} \bm q(\bm{x},t) + \nabla \cdot \bm f (\bm q (\bm{x},t)) &= \bm 0, \label{eq:general_CL}\\
		\bm q(\bm{x},0) &= \bm q_0(\bm{x}),
\end{align}
with initial data $\bm q_0(\bm{x})$ and $\bm{x}\in\Omega$.
We restrict our formulation of Active Flux to a Cartesian grid of $\Omega$ with cells $\mathcal{C}_{i,j,k} = [ x_{i-\frac{1}{2}}, x_{i+\frac{1}{2}} ] \times [ y_{j-\frac{1}{2}}, y_{j+\frac{1}{2}} ] \times [ z_{k-\frac{1}{2}}, z_{k+\frac{1}{2}} ]$, 
$i = 1,\ldots,N_x$, $j = 1,\ldots, N_y$, $k = 1,\ldots, N_z$ and constant mesh widths $\Delta x = x_{i+\frac{1}{2}} - x_{i-\frac{1}{2}}$, $\Delta y = y_{j+\frac{1}{2}} - y_{j-\frac{1}{2}}$ and $\Delta z = z_{k+\frac{1}{2}} - z_{k-\frac{1}{2}}$. 
Integrating (\ref{eq:general_CL}) over one cell $\mathcal{C}_{i,j,k}$, and a timestep $[t_n,t_{n+1}]$ with $\Delta t = t_{n+1} - t_n$, one obtains
\begin{align}
	\bm 0 &= \int_{t_n}^{t_{n+1}} \int_{\mathcal{C}_{i,j,k}} \frac{\partial}{ \partial t} \bm q(\bm{x},t) + \nabla \cdot \bm f (\bm q (\bm{x},t))  \,\dd\bm x \,\dd t \notag \\
	  &= \int_{\mathcal{C}_{i,j,k}} \bm q(\bm{x},t_{n+1}) \,\dd\bm{x} - \int_{\mathcal{C}_{i,j,k}} \bm q(\bm{x},t_{n}) \,\dd\bm{x}  \label{eq:integral_CL}\\
	  &\phantom{=}+ \int_{t_n}^{t_{n+1}} \int_{\mathcal{C}_{i,j,k}} \nabla \cdot \bm f (\bm q (\bm{x},t))  \,\dd\bm x \,\dd t. \notag
\end{align}
Applying the divergence theorem with outer normal vector $\bm{\nu}$ in the latter integral one obtains
\begin{equation}\label{eq:exact_fluxes}
\begin{aligned}
	\int_{\mathcal{C}_{i,j,k}} &\nabla \cdot \bm f (\bm q (\bm{x},t))  \,\dd\bm x = \int_{\partial \mathcal{C}_{i,j,k}} \bm f (\bm q (\bm{x},t)) \cdot \bm{\nu}  \,\dd S_{\bm x}  \\
	&= \phantom{+} \int_{\mathcal{C}_{j,k}^{2D}} f(\bm q(x_{i+\frac{1}{2}},y,z,t)) \,\dd y \,\dd z - \int_{\mathcal{C}_{j,k}^{2D}} f(\bm q(x_{i-\frac{1}{2}},y,z,t)) \,\dd y \,\dd z  \\
	&\phantom{=} + \int_{\mathcal{C}_{i,k}^{2D}} g(\bm q(x,y_{j+\frac{1}{2}},z,t)) \,\dd x \,\dd z - \int_{\mathcal{C}_{i,k}^{2D}} g(\bm q(x,y_{j-\frac{1}{2}},z,t)) \,\dd x \,\dd z \\
	&\phantom{=} + \int_{\mathcal{C}_{i,j}^{2D}} h(\bm q(x,y,z_{k+\frac{1}{2}},t)) \,\dd x \,\dd y - \int_{\mathcal{C}_{i,j}^{2D}} h(\bm q(x,y,z_{k-\frac{1}{2}},t)) \,\dd x \,\dd y, 
\end{aligned}
\end{equation}
with
\begin{align*}
\mathcal{C}_{i,j}^{2D} :&= [x_{i-\frac{1}{2}},x_{i+\frac{1}{2}}] \times [y_{j-\frac{1}{2}},y_{j+\frac{1}{2}}], \\
\mathcal{C}_{j,k}^{2D} :&= [y_{j-\frac{1}{2}},y_{j+\frac{1}{2}}] \times [z_{k-\frac{1}{2}},z_{k+\frac{1}{2}}], \\
\mathcal{C}_{i,k}^{2D} :&= [x_{i-\frac{1}{2}},x_{i+\frac{1}{2}}] \times [z_{k-\frac{1}{2}},z_{k+\frac{1}{2}}].
\end{align*} 
Inserting (\ref{eq:exact_fluxes}) in (\ref{eq:integral_CL}) and with the cell average
\begin{align}
	\overline{Q}_{i,j,k}^{n} \approx \frac{1}{\Delta x \Delta y \Delta z} \int_{\mathcal{C}_{i,j,k}} \bm q(\bm x, t_n) \,\dd\bm{x}
\end{align}
and the numerical flux 
\begin{align}
	F_{i+\frac{1}{2},j,k} \approx \frac{1}{\Delta y \Delta z \Delta t} \int_{t_n}^{t_{n+1}} \int_{\mathcal{C}_{j,k}^{2D}} f(\bm q(x_{i+\frac{1}{2}},y,z,t))
	\,\dd z \,\dd y \,\dd t,
\end{align}
as well as $G_{i,j+\frac{1}{2},k}$ and $H_{i,j,k+\frac{1}{2}}$ analogously, one obtains the update of the fully discrete method
\begin{align*}
	\overline{Q}_{i,j,k}^{n+1} = \overline{Q}_{i,j,k}^{n} & - \frac{\Delta t}{\Delta x} \left(F_{i+\frac{1}{2},j,k} - F_{i-\frac{1}{2},j,k} \right) - 
	\frac{\Delta t}{\Delta y} \left(G_{i,j+\frac{1}{2},k} - G_{i,j-\frac{1}{2},k} \right) 	\\
	&  - \frac{\Delta t}{\Delta z} \left(H_{i,j,k+\frac{1}{2}} - H_{i,j,k-\frac{1}{2}} \right). 
\end{align*}
While the update formula resembles an explicit Euler step, it is a fully discrete scheme: the numerical fluxes $F, G, H$ approximate the flux integrals averaged over $[t_n, t_{n+1}]$
via a quadrature rule, rather than evaluating the flux at a single time level. 
This distinguishes Active Flux from semi-discrete methods where the spatial operator is decoupled from time integration.

\subsection{Degrees of freedom}
\begin{figure}[htbp]
	\centering
	\includegraphics[width = 0.49\textwidth]{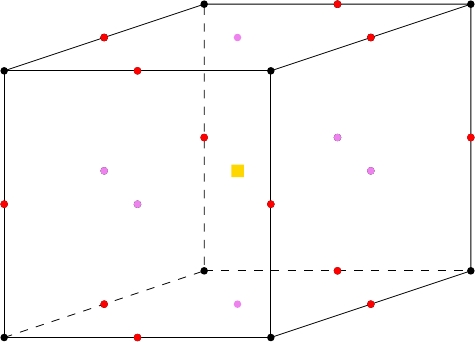}
	\caption{Configuration of degrees of freedom for a three-dimensional cell. 
			 Point values at the nodes (black), edge midpoints (red), face centers (violet) 
			and cell average (yellow).}
	\label{fig:3Dcell}
\end{figure}

Some classical finite volume methods achieve higher-order accuracy through reconstruction from cell averages alone, using wider stencils such as 
WENO \cite{LIU1994200, JIANG1996202} or employing \textit{correction terms} as in wave propagation algorithms \cite{leveque_wave_1997, helzel_third-order_2020}.
Active Flux methods instead introduce point values located on cell boundaries as additional degrees of freedom which are shared between adjacent cells, yielding a globally continuous reconstruction. 
Furthermore, these point values allow for a higher-order method while retaining a compact stencil in space and time.
In our proposed methods we use cell average values, nodal point values and point values on the center of edges and faces. An illustration is given
in Figure~\ref{fig:3Dcell}. There the point values are depicted as circles, while the cell average is presented as a small cube in the center of the cell. 
We denote the nodal point values by
\begin{align}
	Q_{i\pm \frac{1}{2}, j \pm \frac{1}{2}, k \pm \frac{1}{2}}^n \approx \bm q(x_{i \pm \frac{1}{2}}, y_{j \pm \frac{1}{2}}, z_{k \pm \frac{1}{2}}, t_n)
\end{align}
the edge-centered point values by
\begin{align}
	Q_{i , j \pm \frac{1}{2}, k \pm \frac{1}{2}}^n \approx \bm q(x_{i}, y_{j \pm \frac{1}{2}}, z_{k \pm \frac{1}{2}}, t_n),
\end{align}
and $Q_{i \pm \frac{1}{2} , j , k \pm \frac{1}{2}}^n$, $Q_{i \pm \frac{1}{2} , j \pm \frac{1}{2}  , k}^n$, respectively, as well as
the face-centered point values by
\begin{align}
	Q_{i , j \pm \frac{1}{2}, k}^n \approx \bm q(x_{i}, y_{j \pm \frac{1}{2}}, z_{k }, t_n),
\end{align}
and $Q_{i \pm \frac{1}{2} , j , k }^n$, $Q_{i , j , k\pm \frac{1}{2}}^n$, respectively.

\subsection{Numerical fluxes}
In contrast to other finite volume methods, the numerical fluxes are direct approximations of the flux integrals (\ref{eq:exact_fluxes}).
The choice of quadrature rule in space is driven by the aspect that it should use all point values on the boundary. 
The time quadrature should be chosen such that the method is of higher order and not too complex so that point values can be evaluated without excessive computational cost. 
The canonical choice with this configuration of degrees of freedom is Simpson's rule in space and time.
In \cite{AFSchemesEymannRoe}, Roe justified this choice in the one-dimensional case by arguing that it yields a third-order accurate method in space and time. 
Point values at full time steps are used both for the reconstruction and flux computation.
Hence, we will use for the numerical fluxes two-dimensional Simpson's rule in space and one-dimensional Simpson's rule in time as it is the canonical choice for the configuration of the degrees of freedom presented in 
Figure~\ref{fig:3Dcell}.
This gives the following representations of the numerical fluxes:
\begin{equation}
  \begin{aligned}
    F_{i+\frac{1}{2},j,k} = \frac{1}{6} 
    \sum_{(\delta_j,\delta_k)} \omega_{\delta_j,\delta_k} \bigg(
    &f\left(Q_{i+\frac{1}{2}, j+\delta_j, k+\delta_k}^n\right) \\
    + 4 &f\left(Q_{i+\frac{1}{2}, j+\delta_j, k+\delta_k}^{n+\frac{1}{2}}\right)
    + f\left(Q_{i+\frac{1}{2}, j+\delta_j, k+\delta_k}^{n+1}\right) \bigg),
  \end{aligned}
\end{equation}
with $\delta_j, \delta_k \in \{-\frac{1}{2}, 0, \frac{1}{2}\}$,
$\omega_{\delta_j,\delta_k} = \omega_{\delta_j} \omega_{\delta_k}$ and $\omega_{-\frac{1}{2}} = \omega_{\frac{1}{2}} = \frac{1}{6}$, $\omega_{0} = \frac{4}{6}$. The remaining numerical fluxes 
are approximated accordingly.

\subsection{3D reconstruction}
As the formulation of the numerical fluxes suggests, we rely on 26 point values in addition to the cell average. Those point values are shared
between neighboring cells and thus the reconstruction is globally continuous. This stands in contrast to the majority of finite volume and discontinuous Galerkin methods \cite{Hesthaven2008},
which are based on discontinuities between adjacent cells. Consequently, Active Flux methods do not require any form of Riemann solver and the numerical fluxes can be
computed as described above. Roe argued that Riemann solvers introduce one-dimensional physics and do not capture \textit{truly multidimensional} wave behavior \cite{roe_is_2017}. 
Active Flux instead relies on a globally continuous reconstruction and multidimensional evolution operators, which we discuss in the following subsection.
To formulate the three-dimensional reconstruction within one cell, we initially consider the reconstruction over a one-dimensional cell $[x_{i-\frac{1}{2}}, x_{i+\frac{1}{2}}]$. 
With the coordinate transformation 
\begin{align*}
	\xi := \xi(x) = \frac{2(x - x_{i-\frac{1}{2}})}{\Delta x} - 1 \in [-1,1]
\end{align*} 
we use the three Lagrangian basis polynomials
\begin{align*}
	\ell_{\xi, -\frac{1}{2}} = \frac{\xi(\xi - 1)}{2}, \quad \ell_{\xi, 0} =  (1 - \xi^2), \quad \ell_{\xi,\frac{1}{2}} =  \frac{\xi(\xi + 1)}{2}
\end{align*}
that interpolate the left and right boundaries.
Defining $\eta(y)$, $\nu(z)$ and their corresponding Lagrangian polynomials $\ell_{\eta}, \ell_{\nu}$ we obtain, following the tensor-product approach
\begin{align*}
	q^{rec,n}_{i,j,k}(\xi,\eta,\nu) = \sum_{(\hat{i}, \hat{j},\hat{k}) \in \mathcal{I}} Q^n_{i + \hat{i},j + \hat{j}, k + \hat{k}} \ell_{\xi,\hat{i}}\ell_{\eta,\hat{j}} \ell_{\nu,\hat{k}} + \overline{c}_{i,j,k} \ell_{\xi,0}\ell_{\eta,0} \ell_{\nu,0},
\end{align*}
where the index set
\begin{align*}
	\mathcal{I} := \left\{ -\frac{1}{2}, 0, \frac{1}{2} \right\}^3 \setminus \left\{ (0,0,0) \right\}
\end{align*}
enumerates the $26$ point values on the boundary of $\mathcal{C}_{i,j,k}$. The remaining tensor-product basis function $\ell_{\xi,0}\ell_{\eta,0}\ell_{\nu,0}$ is not associated with a point value. Being the only one that 
vanishes on the whole cell boundary, it does not interfere with the interpolation of the point values, and its coefficient $\overline{c}_{i,j,k}$ is instead determined by the cell average, i.e., $\overline{c}_{i,j,k}$ is chosen such that
\begin{align*}
	\overline{Q}_{i,j,k}^n = \frac{1}{2^3}  \int_{-1}^{1} \int_{-1}^{1} \int_{-1}^{1} q^{rec,n}_{i,j,k}(\xi,\eta,\nu) \,\dd\xi \,\dd\eta \,\dd\nu,
\end{align*}
i.e.
\begin{align*}
	\overline{c}_{i,j,k} = \frac{27}{8} \bigg( \overline{Q}_{i,j,k}^n
	- \sum_{(\hat{i},\hat{j},\hat{k}) \in \mathcal{I}} \omega_{\hat{i}}\, \omega_{\hat{j}}\, \omega_{\hat{k}}\,
	Q^n_{i+\hat{i},j+\hat{j},k+\hat{k}} \bigg).
\end{align*}
Thus the used reconstruction fulfills the interpolation on every point value along the boundary of each cube and its average over the cell equals the cell-average. 

\subsection{Evolution operators}
As one can see in the approximation of the numerical fluxes with Simpson's rule, to make a fully discrete Active Flux method work, one has to be able to
give reasonable approximations to the point values at intermediate time level $t_{n+\frac{1}{2}}$ as well as full time level $t_{n+1}$. For the advection equation this is simple characteristic tracing. For more complex PDEs, such a simple approach is not available or diminishes the order of accuracy.
In publications, different types of both exact and approximate evolution operators have been discussed \cite{AFSchemesEymannRoe, eymann_multidimensional_2013, Fan_Roe_2015_wave, chudzik_active_2024, chudzik_fully_2025, barsukow_active_2025}.
So far they have exclusively been used in one- and two-dimensional Active Flux methods. 
Finding a suitable evolution operator is therefore the main difficulty in constructing fully discrete Active Flux methods. 
It should ensure that the higher order is not achieved by information from cells beyond the immediate stencil, and keep the stencil compact \cite{tyacke_musings_2025}. 
Additionally, a multidimensional evolution operator incorporates \textit{truly multidimensional physics} into the method, i.e., 
it accounts for all local propagation directions derived directly from the PDE \cite{chudzik_active_2024}.

\section{An Active Flux method for Maxwell's equations}\label{sec:Maxwell}
\subsection{The mathematical model}
In this article, we consider the source-free Maxwell's equations. They are given by the following system of partial differential equations. 
\begin{subequations}\label{eq: MaxwellVac}
	\begin{align} 
		\nabla \times \bm E &= - \frac{\partial}{\partial t} \bm B, \label{eq:FaradayVac} \\ 
		\nabla \times \bm B &= \mu_0  \varepsilon_0 \frac{\partial}{\partial t} \bm E , \label{eq:AmpereMaxwellVac} \\		
		\nabla \cdot  \bm E &= 0, \label{eq:MagmonVac} \\ 
		\nabla \cdot \bm B &= 0. \label{eq:GaussVac}
	\end{align}
\end{subequations}
Here, $\bm B, \bm E : \mathbb{R}^3 \times \mathbb{R}^+ \to \mathbb{R}^3$  are the magnetic and electric field variables, $\mu_0, \varepsilon_0 \in \mathbb{R}^+$ are the permeability and permittivity satisfying the relation $\mu_0 \varepsilon_0 = c^{-2}$ for the 
speed of light $c$. 
To apply the Active Flux method, we formulate the set of equations (\ref{eq:FaradayVac}--\ref{eq:AmpereMaxwellVac}) in conservative form. Defining the vector 
\begin{align}
	\bm q(\bm x, t) = \begin{pmatrix}
		 B_1 (\bm x, t)\\ 
		 B_2 (\bm x, t)\\
		 B_3 (\bm x, t) \\
		 E_1 (\bm x, t)\\
		 E_2 (\bm x, t)\\
		 E_3 (\bm x, t)
	\end{pmatrix} \in \mathbb{R}^6,
\end{align}
and the matrices 
\begin{align*}
	&\bm A_1 = \begin{pmatrix}
		0 & 0 & 0 & 0 & 0 & 0 \\
		0 & 0 & 0 & 0 & 0 & -1 \\
		0 & 0 & 0 & 0 & 1 & 0 \\
		0 & 0 & 0 & 0 & 0 & 0 \\
		0 & 0 & c^2 & 0 & 0 & 0 \\
		0 & -c^2 & 0 & 0 & 0 & 0 \\
	\end{pmatrix},
	\quad \bm A_2 = \begin{pmatrix}
		0 & 0 & 0 & 0 & 0 & 1 \\
		0 & 0 & 0 & 0 & 0 & 0 \\
		0 & 0 & 0 & -1 & 0 & 0 \\
		0 & 0 & -c^2 & 0 & 0 & 0 \\
		0 & 0 & 0 & 0 & 0 & 0 \\
		c^2 & 0 & 0 & 0 & 0 & 0 \\
	\end{pmatrix}, \\
	& \phantom{--------}\bm A_3 = \begin{pmatrix}
		0 & 0 & 0 & 0 & -1 & 0 \\
		0 & 0 & 0 & 1 & 0 & 0 \\
		0 & 0 & 0 & 0 & 0 & 0 \\
		0 & c^2 & 0 & 0 & 0 & 0 \\
		-c^2 & 0 & 0 & 0 & 0 & 0 \\
		0 & 0 & 0 & 0 & 0 & 0 \\
	\end{pmatrix},
\end{align*}
we can formulate 
\begin{align}\label{eq:MaxwellQuasilinear}
	\partial_t \bm q(\bm x, t) + \bm A_1  \partial_{x} \bm q(\bm x, t) + \bm A_2  \partial_{y} \bm q(\bm x, t) + \bm A_3  \partial_{z} \bm q(\bm x, t) = \bm 0, \\
	\nabla \cdot \bm E(\bm x, t) = 0, \quad \nabla \cdot \bm B(\bm x, t) = 0. \nonumber
\end{align}
Now (\ref{eq:MaxwellQuasilinear}) allows the straightforward application of our Active Flux method.
To give a full description of the method, we have to consider three things as mentioned in the previous section. These are the formulation of the numerical fluxes, the update of the point values, and the reconstruction.
We do not discuss the reconstruction further for this method, as it follows directly from Section~\ref{sec:AF}.
Later, in Section~\ref{sec:divfree}, we will give a description of an alternative reconstruction that is locally divergence-free.
\subsection{Numerical fluxes}
First, we consider the cell-average update, that is, the description of the numerical flux. 
For the six-dimensional flux vector on an arbitrary three-dimensional cell $\mathcal{C}_{i,j,k} = [x_{i-\frac{1}{2}},x_{i+\frac{1}{2}}] \times [y_{j-\frac{1}{2}},y_{j+\frac{1}{2}}]  \times  [z_{k-\frac{1}{2}},z_{k+\frac{1}{2}}] $ over a time step $[t_n,t_{n+1}]$, we have, as described in the previous section, for the face $x_{i+\frac{1}{2}}  \times [y_{j-\frac{1}{2}},y_{j+\frac{1}{2}}]  \times  [z_{k-\frac{1}{2}},z_{k+\frac{1}{2}}]$ of the cell $\mathcal{C}_{i,j,k}$
\begin{equation}
	\begin{aligned}
		&\frac{1}{\Delta y \Delta z \Delta t} \int_{t_n}^{t_{n+1}} \int_{y_{j- \frac{1}{2}}}^{y_{j + \frac{1}{2}}} \int_{z_{k- \frac{1}{2}}}^{z_{k + \frac{1}{2}}} \bm A_1 \bm q (x_{i+\frac{1}{2}}, y, z,t) \,\dd z \,\dd y \,\dd t \\
		& = \frac{1}{\Delta y \Delta z \Delta t} \begin{pmatrix}
			0 \\
			-\int_{t_n}^{t_{n+1}} \int_{y_{j- \frac{1}{2}}}^{y_{j + \frac{1}{2}}} \int_{z_{k- \frac{1}{2}}}^{z_{k + \frac{1}{2}}} E_3 \,\dd z \,\dd y \,\dd t \\
			\int_{t_n}^{t_{n+1}} \int_{y_{j- \frac{1}{2}}}^{y_{j + \frac{1}{2}}} \int_{z_{k- \frac{1}{2}}}^{z_{k + \frac{1}{2}}} E_2 \,\dd z \,\dd y \,\dd t \\
			0 \\
			c^2 \int_{t_n}^{t_{n+1}} \int_{y_{j- \frac{1}{2}}}^{y_{j + \frac{1}{2}}} \int_{z_{k- \frac{1}{2}}}^{z_{k + \frac{1}{2}}}  B_3 \,\dd z \,\dd y \,\dd t \\
			-c^2 \int_{t_n}^{t_{n+1}} \int_{y_{j- \frac{1}{2}}}^{y_{j + \frac{1}{2}}} \int_{z_{k- \frac{1}{2}}}^{z_{k + \frac{1}{2}}} B_2 \,\dd z \,\dd y \,\dd t
		\end{pmatrix}
		&\approx F_{i+\frac{1}{2},j,k}.
	\end{aligned}
\end{equation}
Over the face $[x_{i-\frac{1}{2}}, x_{i+\frac{1}{2}}]\times y_{j+\frac{1}{2}}\times[z_{k-\frac{1}{2}}, z_{k+\frac{1}{2}}]$ we have
\begin{align}
	G_{i,j+\frac{1}{2},k} \approx \frac{1}{\Delta x \Delta z \Delta t} \begin{pmatrix}
		\int_{t_n}^{t_{n+1}} \int_{x_{i- \frac{1}{2}}}^{x_{i + \frac{1}{2}}} \int_{z_{k- \frac{1}{2}}}^{z_{k + \frac{1}{2}}} E_3  \,\dd z \,\dd x \,\dd t \\
		0 \\
		-\int_{t_n}^{t_{n+1}} \int_{x_{i- \frac{1}{2}}}^{x_{i + \frac{1}{2}}} \int_{z_{k- \frac{1}{2}}}^{z_{k + \frac{1}{2}}} E_1  \,\dd z \,\dd x \,\dd t \\
		-c^2 \int_{t_n}^{t_{n+1}} \int_{x_{i- \frac{1}{2}}}^{x_{i + \frac{1}{2}}} \int_{z_{k- \frac{1}{2}}}^{z_{k + \frac{1}{2}}} B_3 \,\dd z \,\dd x \,\dd t \\
		0 \\
		c^2 \int_{t_n}^{t_{n+1}} \int_{x_{i- \frac{1}{2}}}^{x_{i + \frac{1}{2}}} \int_{z_{k- \frac{1}{2}}}^{z_{k + \frac{1}{2}}}  B_1 \,\dd z \,\dd x \,\dd t
	\end{pmatrix},
\end{align}
and over the face $[x_{i-\frac{1}{2}}, x_{i+\frac{1}{2}}]\times [y_{j-\frac{1}{2}},y_{j+\frac{1}{2}}]  \times z_{k+\frac{1}{2}}$
\begin{align}
	H_{i,j,k+\frac{1}{2}} \approx \frac{1}{\Delta x \Delta y \Delta t} \begin{pmatrix}
		-\int_{t_n}^{t_{n+1}} \int_{x_{i- \frac{1}{2}}}^{x_{i + \frac{1}{2}}} \int_{y_{j- \frac{1}{2}}}^{y_{j + \frac{1}{2}}}  E_2  \,\dd y \,\dd x \,\dd t \\
		\int_{t_n}^{t_{n+1}} \int_{x_{i- \frac{1}{2}}}^{x_{i + \frac{1}{2}}} \int_{y_{j- \frac{1}{2}}}^{y_{j + \frac{1}{2}}}  E_1  \,\dd y \,\dd x \,\dd t \\
		0 \\
		c^2 \int_{t_n}^{t_{n+1}} \int_{x_{i- \frac{1}{2}}}^{x_{i + \frac{1}{2}}} \int_{y_{j- \frac{1}{2}}}^{y_{j + \frac{1}{2}}}  B_2 \,\dd y \,\dd x \,\dd t \\
		-c^2 \int_{t_n}^{t_{n+1}} \int_{x_{i- \frac{1}{2}}}^{x_{i + \frac{1}{2}}} \int_{y_{j- \frac{1}{2}}}^{y_{j + \frac{1}{2}}}  B_1 \,\dd y \,\dd x \,\dd t \\
		0
	\end{pmatrix}.
\end{align}
Finally, the numerical fluxes are approximated using a higher-order quadrature rule. We continue with the proposed method in the previous section and use Simpson's rule. For one entry we have
\begin{align}
	&\frac{1}{\Delta y \Delta z \Delta t} \int_{t_n}^{t_{n+1}} \int_{y_{j- \frac{1}{2}}}^{y_{j + \frac{1}{2}}} \int_{z_{k- \frac{1}{2}}}^{z_{k + \frac{1}{2}}}  E_3(x_{i+\frac{1}{2}}, y, z,t) \,\dd z \,\dd y \,\dd t \notag \\
	&\approx \frac{1}{6} \sum_{(\delta_j,\delta_k)} \omega_{\delta_j,\delta_k} \bigg( E_{3,i+\frac{1}{2},j+\delta_j,k+\delta_k}^n
	+ 4  E_{3,i+\frac{1}{2},j+\delta_j,k+\delta_k}^{n+\frac{1}{2}} +  E_{3,i+\frac{1}{2},j+\delta_j,k+\delta_k}^{n+1} \bigg), \notag
\end{align}
with $\delta_j, \delta_k \in \{-\frac{1}{2}, 0, \frac{1}{2}\}$,
$\omega_{\delta_j,\delta_k} = \omega_{\delta_j} \omega_{\delta_k}$ and $\omega_{-\frac{1}{2}} = \omega_{\frac{1}{2}} = \frac{1}{6}$, $\omega_{0} = \frac{4}{6}$. 
What remains is to find an evolution of the point values located on the respective faces of the three-dimensional cube  $[x_{i-\frac{1}{2}},x_{i+\frac{1}{2}}] \times [y_{j-\frac{1}{2}},y_{j+\frac{1}{2}}]  \times  [z_{k-\frac{1}{2}},z_{k+\frac{1}{2}}] $. 

\subsection{Update of point values}
To discuss the direct procedure to update the point value degrees of freedom of this Active Flux method, we first derive an analytical formula.
For the update of the point values, one can exploit the structure of the underlying PDE to update the point values independently \cite{AFSchemesEymannRoe}. 
\subsubsection{Analytical derivation of update formula}
Following the ideas of Eymann and Roe \cite{eymann_multidimensional_2013}, we use the underlying structure of a wave equation for Maxwell's equations.
Applying the \textit{curl} operator to Maxwell's equations, we obtain 
\begin{align*}
	\nabla \times \left( \frac{\partial}{\partial t} \bm B + \nabla \times \bm E \right) = \bm 0 \Leftrightarrow \frac{\partial}{\partial t} \left(\nabla \times \bm B \right) + \nabla \times \left( \nabla \times \bm E \right) = \bm 0, \\
	\nabla \times \left(   \frac{\partial}{\partial t} \frac{1}{c^2} \bm E  - \nabla \times \bm B \right) = \bm 0  \Leftrightarrow \frac{1}{c^2}   \frac{\partial}{\partial t} \left(\nabla \times \bm E \right) - \nabla \times \left( \nabla \times \bm B  \right) = \bm 0.
\end{align*}
Using (\ref{eq:FaradayVac}) and (\ref{eq:AmpereMaxwellVac}), the known identity $ \nabla \times \nabla \times \bm V = \nabla\left( \nabla \cdot \bm V \right) - \Delta  \bm V$ for the vector field $\bm V$ and the divergence constraints (\ref{eq:MagmonVac}), (\ref{eq:GaussVac}), we get
\begin{align}
	&\frac{1}{c^2} \frac{\partial^2}{\partial t^2} \bm E + \nabla \underbrace{ \left(\nabla \cdot \bm E \right)}_{= 0} - \Delta \bm E = \bm 0 \Leftrightarrow \frac{\partial^2}{\partial t^2} \bm E - c^2 \Delta \bm E = \bm 0, \\
	- &\frac{1}{c^2} \frac{\partial^2}{\partial t^2} \bm B - \nabla \underbrace{ \left(\nabla \cdot \bm B \right)}_{= 0} + \Delta \bm B = \bm 0 \Leftrightarrow \frac{\partial^2}{\partial t^2} \bm B - c^2 \Delta \bm B  = \bm 0. 
\end{align}
So we transform Maxwell's equations to the decoupled electromagnetic wave equations. 
To discuss the evolution of our point-values in time, we first recall how to solve a wave equation. Consider the following general wave equation for a function $u$.
\begin{align*}
	\partial_{tt} u(\bm x, t) - c^2\Delta u(\bm x, t) &= 0, \\
	\partial_t u(\bm x, 0) &=g(\bm x), \\ 
	u(\bm x, 0) &=h(\bm x).
\end{align*}
For this set of equations there is a classical solution 
\begin{align}
	u(\bm x, t) = \partial_t \left( t M \{ h, \bm x;ct \}   \right) + t M\{ g,\bm x; ct  \} \label{eq:Kirchhoff}
\end{align}
given by Kirchhoff's formula in three dimensions, see \cite{evans_partial_2010}.
Here, $M\{f,\bm x;r\}$ is the spherical mean of a function $f$ centered at $\bm x$ with radius $r$, i.e.,
\begin{align}
	M\{f,\bm x;r\} := \frac{1}{4\pi r^2} \int_{\partial B_r(\bm x)} f(\bm y) \,\dd S_{\bm{y}}.
\end{align}
Here, $\partial B_r(\bm x)$ is the boundary of the ball of radius $r$ centered at $\bm x$ and $\dd S_{\bm{y}}$ is the area element on the surface of the ball.
In this three-dimensional setting, we use the transformation of the surface integral with spherical coordinates
\begin{align*}
	M\{f,\bm x;r\} = \frac{1}{4\pi } \int_{0}^{2 \pi} \int_{0}^{\pi}  f(& x +r\sin(\phi)\cos(\theta) , y +r\sin(\phi)\sin(\theta) , \\
	& z+r\cos(\phi)) \sin(\phi	) \,\dd\phi \,\dd\theta.
\end{align*}
In total we can formulate the following two sets of wave equations for the magnetic field
\begin{equation}\label{eq:EM_WaveB}
	\begin{aligned}
		\frac{\partial^2}{\partial t^2} \bm B - c^2 \Delta \bm B &=  \bm 0, \\
		\frac{\partial}{\partial t} \bm B (\bm x, 0) &=  - \nabla \times \bm E (\bm x, 0)  , \\
		\bm B (\bm x, 0) &= \bm B_0(\bm x),
	\end{aligned}
\end{equation}
and for the electric field,
\begin{equation}
	\begin{aligned}
		\frac{\partial^2}{\partial t^2} \bm E - c^2 \Delta \bm E &= \bm 0, \\ 
		\frac{\partial}{\partial t} \bm E (\bm x, 0)   &=   c^2  \nabla \times \bm B (\bm x, 0), \\
		\bm E (\bm x, 0) &= \bm E_0(\bm x). 
	\end{aligned}
\end{equation}
For sufficiently smooth solutions of the wave equations with divergence-free initial data, one can show $\partial_t \nabla \cdot \bm B = 0$ and $\partial_t \nabla \cdot \bm E = 0$.
Therefore, numerical methods based on exact evolution operators for the wave equation are interesting candidates for structure--preserving methods. 
Following the wave equation ansatz from above, we can derive with (\ref{eq:Kirchhoff}) and (\ref{eq:EM_WaveB}) the update formula for the components $B_m$, $m = 1,2,3$,
\begin{align*}
	 B_m(\bm x, t) &= \partial_t \left( t M\{ (\bm B_0)_m,\bm x;ct   \}   \right) + t  M\{ - (\nabla \times \bm E_0 )_m,\bm x;ct   \}    \\
	&=  M \{  (\bm B_0)_m, \bm{x};ct  \} + t \partial_t  M \{   (\bm B_0)_m , \bm{x};ct \} - t M\{   ( \nabla \times \bm E_0 )_m, \bm{x};ct \}. 
\end{align*}
Defining $ct =: R$ and following the calculation from Eymann and Roe in \cite{eymann_multidimensional_2013}, we can rewrite the second term
by a change of variables and obtain
\begin{align*}
	t \partial_t  M \{ (\bm B_0)_m,\bm{x};R  \} = t  \partial_R  M \{ (\bm B_0)_m , \bm{x};R \} \frac{\partial R}{\partial  t} = R \partial_R  M \{ (\bm B_0)_m , \bm{x};R  \}.
\end{align*}
With these calculations, we get the evolution formula 
\begin{equation}
	\begin{aligned}
		B_m(\bm x, t)  &=   M \{ (\bm B_0)_m ,\bm{x};R  \}  +  R \partial_R  M \{  (\bm B_0)_m , \bm{x};R \} \\
		&\quad - t M \{ ( \nabla \times \bm E_0 )_m , \bm{x};R\}. \label{eq:continuous_updateB}
	\end{aligned}
\end{equation}
Following the same arguments, one derives the formula for the electric field 
\begin{equation}
	\begin{aligned}
		 E_m(\bm x, t)  &=   M \{ (\bm E_0)_m  , \bm{x};R\} + R \partial_R  M \{  (\bm E_0)_m ,\bm{x};R  \}  \\
		&\quad+ c^2 t M \{ ( \nabla \times \bm B_0  )_m, \bm{x};R \}. \label{eq:continuous_updateE}
	\end{aligned}
\end{equation}

\subsubsection{Discrete update of point values}
Using the derived evolution formula, we demonstrate the update of a point value
located at $\bm x$ on the boundary of the cell $\mathcal{C}_{i,j,k}$.
Using the previously derived formulas with $R = c \Delta t$, we proceed as follows:
\begin{align*}
	B_{m;i \pm \frac{1}{2}, j \pm  \frac{1}{2}, k\pm  \frac{1}{2}}^{n+1} = &\phantom{-} M \{ B^n_{m},  \bm x ; R   \}
	+ R \partial_R  M \{ B^n_{m},  \bm x ; R   \}  \\
	&- \Delta t M\{ ( \nabla \times \bm E^n)_m, \bm x ; R   \} ,
	\end{align*}
	and
	\begin{align*}
	E_{m;i \pm \frac{1}{2}, j \pm  \frac{1}{2}, k\pm  \frac{1}{2}}^{n+1} =  &\phantom{-} M \{ E^n_{m},  \bm x ; R   \}
	+ R \partial_R  M \{ E^n_{m},  \bm x ; R   \}  \\
	&+ c^2 \Delta t M\{ ( \nabla \times \bm B^n)_m, \bm x ; R   \},
\end{align*}
for $m = 1,2,3$.
Here, $B_{m; i \pm \frac{1}{2}, j \pm  \frac{1}{2}, k\pm  \frac{1}{2}}^{n+1} \approx B_m(\bm x, t_{n+1})$ and $E_{m;i \pm \frac{1}{2}, j \pm  \frac{1}{2}, k\pm  \frac{1}{2}}^{n+1} \approx E_m(\bm x, t_{n+1})$, 
while $B_m^{n}$ and $E_m^n$ denote the globally continuous reconstructions of the $m$-th component of the magnetic and electric fields for $m=1,2,3$, and $\bm B^{n},  \bm E^{n}$ are the vectors of reconstructions.
Now, to update a point-value used in our method, we need to consider the three different cases from above. First, the general calculation of a spherical mean of one component of a field, second, the derivative of 
the spherical mean with respect to the radius $R$, and third, the spherical mean of a partial derivative of a field component.

\subsubsection{Approximation of spherical means}
For simplicity, we demonstrate the calculations for an arbitrary function $q$ and its associated reconstruction $q^{rec}$ at time 
$t_n$, interpolating the $26$ point values $Q^n_{i+\hat{i}, j+\hat{j}, k+\hat{k}}$, $(\hat{i},\hat{j},\hat{k}) \in \mathcal{I}$, on the boundary and 
conserving the cell average $\overline{Q}_{i,j,k}^n$ in $\mathcal{C}_{i,j,k}$ as described in 
Section~\ref{sec:AF}.
For the general computation of a spherical mean of a piecewise-defined reconstruction, we decompose the integral over the 
entire sphere into integrals over spherical wedges on cells on which the reconstruction is locally defined. Thus, the spherical mean 
of $q$ is the sum of integrals over wedges $W^{\bar{i},\bar{j}, \bar{k}}$
\begin{align*}
	M\{q, \bm x;R\} = \sum_{\bar{i},\bar{j}, \bar{k}} M^{W^{\bar{i},\bar{j}, \bar{k}}}\{q, \bm x;R\},
\end{align*}
where, in this three-dimensional setting, we associate the spherical-wedge $W^{\bar{i},\bar{j}, \bar{k}}$ to $(\phi_{min}^{\bar{i},\bar{j}, \bar{k}},\phi_{max}^{\bar{i},\bar{j}, \bar{k}})\times(\theta_{min}^{\bar{i},\bar{j}, \bar{k}},\theta_{max}^{\bar{i},\bar{j}, \bar{k}})$ and
\begin{align*}
	M^{W^{\bar{i},\bar{j}, \bar{k}}}\{q,\bm x;R\} := \frac{1}{4\pi} \int_{\theta_{min}^{\bar{i},\bar{j}, \bar{k}}}^{\theta_{max}^{\bar{i},\bar{j}, \bar{k}}} \int_{\phi_{min}^{\bar{i},\bar{j}, \bar{k}}}^{\phi_{max}^{\bar{i},\bar{j}, \bar{k}}} q(& x +R\sin(\phi)\cos(\theta) , \\
	& y + R\sin(\phi)\sin(\theta) , \\
	&z  + R\cos(\phi)) \sin(\phi) \,\dd\phi \,\dd\theta.
\end{align*}
When it comes to integrating the reconstruction $q^{rec}$, we first consider the coordinate transformations arising from the arguments above
\begin{align*}
	\xi &= \xi(\phi, \theta, R, x) = 2\frac{( x + R\sin(\phi)\cos(\theta) - x_{\bar{i}-\frac{1}{2}})}{\Delta x} - 1, \\
	\eta &= \eta(\phi, \theta, R, y) = 2\frac{( y + R\sin(\phi)\sin(\theta) - y_{\bar{j}-\frac{1}{2}})}{\Delta y} - 1,  \\
	\nu &= \nu(\phi, R, z) = 2\frac{( z + R\cos(\phi) - z_{\bar{k}-\frac{1}{2}})}{\Delta z} - 1,
\end{align*}
that entail the dependencies on $\phi$, $\theta$ and $R$.
The spherical mean of $q^{rec}$ centered at the point $(x,y,z)$, coinciding with the position of a point value on the boundary,
and radius $R$, is given by
\begin{equation}
	\begin{aligned}
		M\{q^{rec,n}, (x,y,z); R\} &= \sum_{\bar{i},\bar{j},\bar{k}} M^{W^{\bar{i},\bar{j},\bar{k}}}\{q^{rec,n}, (x,y,z);R\} \\
		&= \frac{1}{4\pi} \sum_{\bar{i},\bar{j},\bar{k}} \int_{W^{\bar{i},\bar{j},\bar{k}}}
		\bigg( \sum_{(\hat{i},\hat{j},\hat{k}) \in \mathcal{I}} Q^n_{\bar{i}+\hat{i},\bar{j}+\hat{j},\bar{k}+\hat{k}} \, \ell_{\xi,\hat{i}}\ell_{\eta,\hat{j}}\ell_{\nu,\hat{k}} \\
		&\qquad + \overline{c}_{\bar{i},\bar{j},\bar{k}} \, \ell_{\xi,0}\ell_{\eta,0}\ell_{\nu,0} \bigg) \sin(\phi) \,\dd\phi \,\dd\theta \\
		&= \frac{1}{4\pi} \sum_{\bar{i},\bar{j},\bar{k}} \bigg( \sum_{(\hat{i},\hat{j},\hat{k}) \in \mathcal{I}} Q^n_{\bar{i}+\hat{i},\bar{j}+\hat{j},\bar{k}+\hat{k}}
		\int_{W^{\bar{i},\bar{j},\bar{k}}} \ell_{\xi,\hat{i}}\ell_{\eta,\hat{j}}\ell_{\nu,\hat{k}} \sin(\phi) \,\dd\phi \,\dd\theta \\
		&\qquad + \overline{c}_{\bar{i},\bar{j},\bar{k}} \int_{W^{\bar{i},\bar{j},\bar{k}}} \ell_{\xi,0}\ell_{\eta,0}\ell_{\nu,0} \sin(\phi) \,\dd\phi \,\dd\theta \bigg) \\
		&= \frac{1}{4\pi} \bm{d}^n(x,y,z)^T \bm{\ell}(R). \label{eq:SpherMeanqrec}
	\end{aligned}
\end{equation}
where $\bm{d}^n(x,y,z)$ denotes the vector of all degrees of freedom associated with the point $(x,y,z)$, i.e., the point values and the coefficients of all adjacent cells, 
and $\bm{\ell}(R)$ the vector of the associated integrated basis functions, which depends only on $R = c \Delta t$.
The integration reduces to integrating trigonometric polynomials, which we carried out beforehand analytically using SymPy.
Thus, updating the point values requires recomputing $\bm{\ell}(R)$ and evaluating this product.

For the derivatives $\partial_R  M \{q^{rec, n}, \bm x; R \} $ we have with the notation from above 
\begin{align}
\partial_R  M \{q^{rec, n}, \bm x; R \} = \frac{1}{4 \pi} \bm{d}^n(x,y,z)^T \partial_R \bm{\ell}(R), 
\end{align}
where the differentiation is component-wise. The vector $\partial_R \bm{\ell}(R)$ is computed analytically beforehand using SymPy as well.

For the differentiation with respect to $x,y,z$ one must apply the chain rule. We have
\begin{align*}
	\partial_x ( \ell_{\xi,\hat{i}}\ell_{\eta,\hat{j}} \ell_{\nu,\hat{k}} ) = \ell'_{\xi,\hat{i}} \xi'(x) \ell_{\eta,\hat{j}} \ell_{\nu,\hat{k}} = \ell'_{\xi,\hat{i}} \frac{2}{\Delta x} \ell_{\eta,\hat{j}}  \ell_{\nu,\hat{k}}.
\end{align*}
Now we plug this polynomial into the expression (\ref{eq:SpherMeanqrec}). Analogously for the other two derivatives $\partial_y (\ell_{\xi,\hat{i}}\ell_{\eta,\hat{j}} \ell_{\nu,\hat{k}})$ and $\partial_z (\ell_{\xi,\hat{i}}\ell_{\eta,\hat{j}} \ell_{\nu,\hat{k}})$. 
These give us the last terms in the evolution formulas (\ref{eq:continuous_updateB}) and (\ref{eq:continuous_updateE}).

A natural bound of the time step $\Delta t$ by the CFL condition is given by the fact that we only want to integrate over the cells adjacent to the point value.
This means that CFL $\leq \frac{1}{2}$ is necessary to keep the stencil compact, i.e., $R = c \Delta t \leq \frac{1}{2} \min(\Delta x, \Delta y, \Delta z)$.
For two-dimensional Active Flux methods for the acoustic equations this condition is also sufficient for stability \cite{chudzik_active_2024}.
Our numerical results indicate that the same holds for the three-dimensional method presented here.

One time step of the method is summarized in Algorithm~\ref{alg:classical}.

\begin{algorithm}[htbp]
\caption{Active Flux method of Section~\ref{sec:Maxwell}.}
\label{alg:classical}
\begin{algorithmic}[1]
\State \textbf{Reconstruct} $\bm B^{rec,n}_{i,j,k}$ and $\bm E^{rec,n}_{i,j,k}$ on every cell (Section~\ref{sec:AF}).
\State \textbf{Evolve point values:} $\bm{B}^{n+\frac{1}{2}}$, $\bm{E}^{n+\frac{1}{2}}$, $\bm{B}^{n+1}$ and $\bm{E}^{n+1}$ from
  (\ref{eq:continuous_updateB}), (\ref{eq:continuous_updateE}) with $R = \frac{1}{2} c \Delta t$ and $R = c \Delta t$.
\State \textbf{Numerical fluxes:} Compute $F_{i+\frac{1}{2},j,k}$, $G_{i,j+\frac{1}{2},k}$, $H_{i,j,k+\frac{1}{2}}$ with Simpson's rule
  in space and time from the face point values at $t_n$, $t_{n+\frac{1}{2}}$ and $t_{n+1}$.
\State \textbf{Update cell averages:} With $\overline{\bm{Q}}^n_{i,j,k} = (\overline{\bm{B}}^{n}_{i,j,k}, \overline{\bm{E}}^{n}_{i,j,k})^T$ 
  \Statex \hspace{\algorithmicindent}$\begin{aligned}
  \overline{\bm{Q}}^{n+1}_{i,j,k} = \overline{\bm{Q}}^{n}_{i,j,k}
  &- \frac{\Delta t}{\Delta x}\bigl(F_{i+\frac{1}{2},j,k} - F_{i-\frac{1}{2},j,k}\bigr)
  - \frac{\Delta t}{\Delta y}\bigl(G_{i,j+\frac{1}{2},k} - G_{i,j-\frac{1}{2},k}\bigr) \\
  &- \frac{\Delta t}{\Delta z}\bigl(H_{i,j,k+\frac{1}{2}} - H_{i,j,k-\frac{1}{2}}\bigr).
  \end{aligned}$
\end{algorithmic}
\end{algorithm}

\section{A divergence-free method based on the Active Flux approach}\label{sec:divfree}
The method described in the previous section is an extension of the previous Active Flux methods in 2D for acoustics 
or the Euler equations to Maxwell's equations. Although the continuous evolution operators
(\ref{eq:continuous_updateB}), (\ref{eq:continuous_updateE}) preserve $\nabla \cdot \bm B = 0$ and
$\nabla \cdot \bm E = 0$ for sufficiently smooth initial data, the reconstruction $q^{rec,n}$ used in that method
does not satisfy a discrete divergence-free property.
In the following we present an alternative method for Maxwell's equations that is 
divergence-free by construction and uses fewer degrees of freedom. 
For this approach we follow the ideas of \cite{BALSARA2004,BALSARA2009} for divergence-free reconstructions in finite volume schemes for the ideal MHD equations. While Balsara proposed his approaches for MHD, they have been applied to Maxwell's equations as well \cite{HAZRA2019298, BALSARA20191108}. Note that in this source-free regime 
we apply the same ideas to the electric field as well as to the magnetic field. 
Besides the divergence-free property, one of the main differences between the method from Section~\ref{sec:Maxwell} and the one in this section is the reconstruction and its degrees of freedom. The new method uses degrees of freedom 
that are defined on the faces of the cells. We still rely on point values on the cell boundary to evolve the solution in time, but we no longer use them as degrees of freedom for the
reconstruction. 
This differs from Active Flux methods, which use point values on the cell boundary, possibly together with their derivatives, and cell averages as degrees of freedom.

\subsection{Divergence-free reconstruction}\label{sec:divfree_rec}
The building blocks of the divergence-free reconstruction are the degrees of freedom. Instead of placing the degrees of freedom as nodes on the boundary of 
each cell, we use degrees of freedom that are defined on the faces of the cells for each component. 
On a cell we consider the following faces
\begin{align}\label{eq:faces}
	S^{x}_{i-\frac{1}{2},j,k} &:= \{x_{i-\frac{1}{2}}\} \times [y_{j-\frac{1}{2}},y_{j+\frac{1}{2}}] \times [z_{k-\frac{1}{2}},z_{k+\frac{1}{2}}], \nonumber \\
	S^{y}_{i,j-\frac{1}{2},k} &:= [x_{i-\frac{1}{2}},x_{i+\frac{1}{2}}] \times \{y_{j-\frac{1}{2}}\} \times [z_{k-\frac{1}{2}},z_{k+\frac{1}{2}}], \\
	S^{z}_{i,j,k-\frac{1}{2}} &:= [x_{i-\frac{1}{2}},x_{i+\frac{1}{2}}] \times [y_{j-\frac{1}{2}},y_{j+\frac{1}{2}}] \times \{z_{k-\frac{1}{2}}\}. \nonumber
\end{align}
As basis polynomials we use the same ones as in~\cite{BALSARA2009}, rescaled to the interval
$[-1,1]$. Those are given by the following Legendre polynomials
\begin{align}\label{eq:legendre}
	P_0(\xi) = 1, \quad P_1(\xi) = \xi, \quad P_2(\xi) = \xi^2 - \frac{1}{3}, \quad P_3(\xi) = \xi^3 - \frac{3}{5}\xi.
\end{align}
For the face $S^{x}_{i-\frac{1}{2},j,k}$ we use the following ansatz in the reference space $[-1,1]^2$ to describe $B_1$
\begin{equation}\label{eq:face_poly}
\begin{aligned}
	B_1(-1, \eta,\nu) &= \mathcal{B}^{x,-}_{0} P_0(\eta) P_0(\nu) + \mathcal{B}^{x,-}_{y} P_1(\eta) P_0(\nu) + \mathcal{B}^{x,-}_{z} P_0(\eta) P_1(\nu) \\
	&\phantom{=} + \mathcal{B}^{x,-}_{yy} P_2(\eta) P_0(\nu) + \mathcal{B}^{x,-}_{yz} P_1(\eta) P_1(\nu) + \mathcal{B}^{x,-}_{zz} P_0(\eta) P_2(\nu).
\end{aligned}
\end{equation}
This approach is applied analogously to the face $S^{x}_{i+\frac{1}{2},j,k}$, the faces $S^{y}_{i,j-\frac{1}{2},k}$ and $S^{y}_{i,j+\frac{1}{2},k}$ for $B_2$ and 
$S^{z}_{i,j,k-\frac{1}{2}}$ and $S^{z}_{i,j,k+\frac{1}{2}}$ for $B_3$.
In our approach, we approximate these coefficients using updated point values on these faces. 
For each component of the magnetic or electric field, we carry six degrees of freedom. 
On the face $S^{x}_{i-\frac{1}{2},j,k}$, we approximate the coefficients of (\ref{eq:face_poly}) from the point values $B_1(-1,p,q)$, $p,q \in \{-1,0,1\}$, by
\begin{equation}\label{eq:facial_FD}
	\begin{aligned}
		\mathcal{B}^{x,-}_{0}  &\approx \frac{1}{4} \int_{-1}^{1} \int_{-1}^{1} B_1(-1,\eta,\nu) \,\dd\eta \,\dd\nu, \\
		\mathcal{B}^{x,-}_{y}  &:= \frac{1}{2}\big(B_1(-1,1,0) - B_1(-1,-1,0)\big) \approx \partial_{\eta} B_1 \big|_{(-1,0,0)},\\
		\mathcal{B}^{x,-}_{z}  &:= \frac{1}{2}\big(B_1(-1,0,1) - B_1(-1,0,-1)\big) \approx \partial_{\nu} B_1 \big|_{(-1,0,0)},\\
		\mathcal{B}^{x,-}_{yy} &:= \frac{1}{2}\big(B_1(-1,1,0) - 2B_1(-1,0,0) + B_1(-1,-1,0)\big) \approx \frac{1}{2} \partial_{\eta}^2 B_1 \big|_{(-1,0,0)},\\
		\mathcal{B}^{x,-}_{zz} &:= \frac{1}{2}\big(B_1(-1,0,1) - 2B_1(-1,0,0) + B_1(-1,0,-1)\big) \approx \frac{1}{2} \partial_{\nu}^2 B_1 \big|_{(-1,0,0)},\\
		\mathcal{B}^{x,-}_{yz} &:= \frac{1}{4}\big(B_1(-1,1,1) - B_1(-1,1,-1) - B_1(-1,-1,1) + B_1(-1,-1,-1)\big)\\
		&\phantom{:}\approx \partial_{\eta} \partial_{\nu} B_1 \big|_{(-1,0,0)}.
	\end{aligned}
\end{equation}
Analogously for the coefficients on the face $S^{y}_{i,j-\frac{1}{2},k}$ for $B_2$ and $S^{z}_{i,j,k-\frac{1}{2}}$ for $B_3$.
Note that the face average $\mathcal{B}^{x,-}_{0}$ is usually initialized by a vector potential in order to preserve the divergence-free property, however in some cases it suffices to initialize it by Simpson's rule from the point values on the face.
The remaining coefficients are computed from the point values on the face, which are evolved in time by the evolution operator (\ref{eq:continuous_updateB}) with the curl taken in the distributional sense, see (\ref{eq:spherical_jump}).
Besides the fact that this approach leads to a divergence-free reconstruction, it also reduces the number of degrees of freedom. 
For a field we now have to store $3$ faces with $6$ degrees of freedom each, since each face is shared by two cells, giving a total of $18$ facial degrees of freedom per field. 
Counted in the same way, the method of Section~\ref{sec:Maxwell} stores $24$ degrees of freedom for each field per cell.
Given the facial reconstruction, we use the following ansatz for the reconstruction of a component of the magnetic field inside a cell
\begin{equation}\label{eq:divfree_ansatz}
	\begin{aligned}
		B_1^{rec}(\xi, \eta, \nu) = \ & a_0 + a_x \xi + a_y \eta + a_z \nu
		+ a_{xx} P_2(\xi) + a_{yy} P_2(\eta) + a_{zz} P_2(\nu) \\
		& + a_{xy} \xi\eta + a_{xz} \xi\nu + a_{yz} \eta\nu
		+ a_{xyz} \xi\eta\nu \\
		& + a_{xyy} \xi P_2(\eta) + a_{xzz} \xi P_2(\nu)
		+ a_{xxy} P_2(\xi) \eta \\
		& + a_{xxz} P_2(\xi) \nu
		+ a_{xxx} P_3(\xi),
	\end{aligned}
\end{equation}
for $\xi, \eta, \nu \in [-1,1]$. The reconstructions of $B_2^{rec}$ and $B_3^{rec}$ are given in Appendix~\ref{app:reconstruction}, one obtains them by the cyclic permutations.
By matching the reconstruction $B_1^{rec}$, $B_2^{rec}$ and $B_3^{rec}$ to their facial reconstructions as given above in (\ref{eq:face_poly}) and requiring that the divergence is constant in each cell,
\begin{align}\label{eq:pointwise_divfree}
	\partial_\xi B_1^{rec} + \partial_\eta B_2^{rec} + \partial_\nu B_3^{rec} \equiv \text{const}, 
\end{align}
one obtains the desired reconstruction. The value of this constant is not specified. It follows from the reconstruction and is computed 
in Section~\ref{sec:divfree_stokes}. For completeness, we list the resulting coefficients for our rescaling onto $[-1,1]^3$ in Appendix~\ref{app:coefficients}. 

\subsection{Update of the facial degrees of freedom via Stokes' theorem}\label{sec:divfree_stokes}
Computing the divergence of the reconstruction (\ref{eq:divfree_ansatz}) yields
\begin{align}\label{eq:discrete_gauss}
	\nabla_{(\xi, \eta, \nu)} \cdot \bm{B}^{rec} = \frac{1}{2} \left( \mathcal{B}^{x,+}_{0} - \mathcal{B}^{x,-}_{0}
	+ \mathcal{B}^{y,+}_{0} - \mathcal{B}^{y,-}_{0}
	+ \mathcal{B}^{z,+}_{0} - \mathcal{B}^{z,-}_{0} \right)
\end{align}
in each cell.
This expression is a discrete version of Gauss' law for both the magnetic and electric fields. Hence the divergence of the reconstruction vanishes if and only if the right-hand side of 
(\ref{eq:discrete_gauss}) vanishes, i.e. if the face averages satisfy the discrete Gauss law. Consequently, controlling the divergence of the numerical solution reduces entirely to controlling (\ref{eq:discrete_gauss}). 
To ensure that (\ref{eq:discrete_gauss}) remains zero for all times, one needs to compute 
the facial averages obeying Stokes' theorem.
In the following we will describe how we update the facial averages such that the divergence-free property is preserved. 
Multiplying (\ref{eq:FaradayVac}) with the normal vector $\bm{n}_x$ and integrating  over the face $S^{x}_{i-\frac{1}{2},j,k}$ yields
\begin{equation}\label{eq:stokes_exact}
	\begin{aligned}
		&\frac{\dd}{\dd t} \int_{S^{x}_{i-\frac{1}{2},j,k}} \bm B \cdot \bm{n}_x \,\dd S = - \int_{S^{x}_{i-\frac{1}{2},j,k}} \left( \nabla \times \bm E \right) \cdot \bm{n}_x \,\dd S \\
		&\phantom{\frac{\dd}{\dd t} \int_{S^{x}_{i-\frac{1}{2},j,k}} }= - \int_{y_{j-\frac{1}{2}}}^{y_{j+\frac{1}{2}}} \int_{z_{k-\frac{1}{2}}}^{z_{k+\frac{1}{2}}} \left( \partial_y E_3 - \partial_z E_2 \right) \,\dd z \,\dd y \\
		&\phantom{\frac{\dd}{\dd t} \int_{S^{x}_{i-\frac{1}{2},j,k}} }= - \underbrace{\int_{z_{k-\frac{1}{2}}}^{z_{k+\frac{1}{2}}} E_3(x_{i-\frac{1}{2}},y_{j+\frac{1}{2}},z,t) \,\dd z }_{:= C^3_{i-\frac{1}{2}, j+\frac{1}{2}}(t)} + \underbrace{\int_{z_{k-\frac{1}{2}}}^{z_{k+\frac{1}{2}}} E_3(x_{i-\frac{1}{2}},y_{j-\frac{1}{2}},z,t)  \,\dd z}_{:= C^3_{i-\frac{1}{2}, j-\frac{1}{2}}(t)} \\
		&\phantom{\frac{\dd}{\dd t} \int_{S^{x}_{i-\frac{1}{2},j,k}} =} + \underbrace{\int_{y_{j-\frac{1}{2}}}^{y_{j+\frac{1}{2}}} E_2(x_{i-\frac{1}{2}},y, z_{k+\frac{1}{2}},t) \,\dd y}_{:= C^2_{i-\frac{1}{2}, k+\frac{1}{2}}(t)} - \underbrace{\int_{y_{j-\frac{1}{2}}}^{y_{j+\frac{1}{2}}} E_2(x_{i-\frac{1}{2}},y,z_{k-\frac{1}{2}},t)  \,\dd y}_{:= C^2_{i-\frac{1}{2}, k-\frac{1}{2}}(t)}.
	\end{aligned}
\end{equation}
Integrating (\ref{eq:stokes_exact}) in time over a time-step $[t_n,t_{n+1}]$ and dividing by the corresponding area $\Delta y \Delta z$ of the face gives the update formula for 
the facial average $\mathcal{B}^{x,-}_{0}$
\begin{equation}\label{eq:stokes_time}
	\begin{aligned}
		\mathcal{B}^{x,-,n+1}_{0}
		= \mathcal{B}^{x,-,n}_{0}
		- \frac{1}{\Delta y\Delta z}
		\int_{t_n}^{t_{n+1}} \bigg( &C^3_{i-\frac{1}{2}, j+\frac{1}{2}}(t) - C^3_{i-\frac{1}{2}, j-\frac{1}{2}}(t) \\
								- &C^2_{i-\frac{1}{2}, k+\frac{1}{2}}(t) + C^2_{i-\frac{1}{2}, k-\frac{1}{2}}(t) \bigg) \,\dd t.
	\end{aligned}
\end{equation}
So the update of the facial averages is given by the circulation of the electric field along the boundary of the face.
With 
\begin{align}
	C^1_{j-\frac{1}{2}, k-\frac{1}{2}}(t) &:= \int_{x_{i-\frac{1}{2}}}^{x_{i+\frac{1}{2}}} E_1(x,y_{j-\frac{1}{2}}, z_{k-\frac{1}{2}},t) \,\dd x
\end{align}
the remaining formulae for the facial averages $\mathcal{B}^{y,-,n+1}_{0}$ and $\mathcal{B}^{z,-,n+1}_{0}$ are given by
\begin{equation}\label{eq:stokes_time_yz}
	\begin{aligned}
		\mathcal{B}^{y,-,n+1}_{0} = \mathcal{B}^{y,-,n}_{0}
		- \frac{1}{\Delta x\Delta z}\int_{t_n}^{t_{n+1}} \bigg( &C^1_{j-\frac{1}{2}, k+\frac{1}{2}}(t) - C^1_{j-\frac{1}{2}, k-\frac{1}{2}}(t) \\
															- &C^3_{i+\frac{1}{2}, j-\frac{1}{2}}(t) + C^3_{i-\frac{1}{2}, j-\frac{1}{2}}(t) \bigg) \,\dd t, \\
		\mathcal{B}^{z,-,n+1}_{0}
		= \mathcal{B}^{z,-,n}_{0}
		- \frac{1}{\Delta x\Delta y}
		\int_{t_n}^{t_{n+1}} \bigg( &C^2_{i+\frac{1}{2}, k-\frac{1}{2}}(t) - C^2_{i-\frac{1}{2}, k-\frac{1}{2}}(t) \\
								- &C^1_{j+\frac{1}{2}, k-\frac{1}{2}}(t) + C^1_{j-\frac{1}{2}, k-\frac{1}{2}}(t) \bigg) \,\dd t. 
	\end{aligned}
\end{equation}
In order to maintain third-order accuracy, we approximate the line integrals in space with Simpson's rule. For this we use 
point values on the edges of the faces. 
\begin{equation}\label{eq:edge_simpson}
	\begin{aligned}
		C^1_{j-\frac{1}{2}, k-\frac{1}{2}}(t) \approx \frac{\Delta x}{6} \Big(&E_1( x_{i-\frac{1}{2}},y_{j-\frac{1}{2}},z_{k-\frac{1}{2}},t) \\
		+ 4 &E_1(x_{i},y_{j-\frac{1}{2}},z_{k-\frac{1}{2}},t) + E_1( x_{i+\frac{1}{2}},y_{j-\frac{1}{2}},z_{k-\frac{1}{2}},t)\Big).
	\end{aligned}
\end{equation}
The time integration in (\ref{eq:stokes_time}) uses the quadrature rule described in Section~\ref{sec:divfree_time}. The respective point values
on the edges are updated by the evolution operator (\ref{eq:continuous_updateB}) and  (\ref{eq:continuous_updateE}) with distributional curl (\ref{eq:spherical_jump}) at the stages of the time quadrature.
These updates (\ref{eq:stokes_time}), (\ref{eq:stokes_time_yz}) give a divergence-free evolution of the facial averages of $\bm B$ in the sense that if (\ref{eq:discrete_gauss}) 
is satisfied at time $t_n$, it is also satisfied at time $t_{n+1}$. 
Essentially, this update is a discrete version of $\nabla \cdot \left( \nabla \times \bm{E} \right) = 0$, which is applied in many 
structure-preserving methods both for MHD \cite{BALSARA2004, BALSARA2009, FAMBRI2025113691} and Maxwell's equations \cite{BALSARA20191108, KORMANN2025118290, bach2025sbpfdecsummationbypartsfinitedifference}. 

In our setting of considering Maxwell's equations in vacuum, we treat the electric field in the same way as the magnetic field, and the facial averages of $\bm E$ can be updated in
the same way, leading to $\nabla \cdot \bm{E} = 0$.

\subsection{Evolution of point values}
We intend to update the point values used in (\ref{eq:facial_FD}) in the same way as in Section~\ref{sec:Maxwell}, i.e. as described in (\ref{eq:continuous_updateB}), (\ref{eq:continuous_updateE}).
However, this approach is not directly applicable to the reconstruction in this section. 
The reconstruction of Section~\ref{sec:Maxwell} is globally continuous. Thus, computing the 
spherical mean of the curl of $\bm B$ and $\bm E$ is defined in the classical sense. 
The reconstruction (\ref{eq:divfree_ansatz}) is not continuous across grid cell interfaces. 
This has the consequence, that the derivative across grid cell interfaces is not defined in the classical sense. Therefore, we instead consider the distributional derivative over discontinuities normal to the derivative. 
For the spherical mean of the field variables as well as the derivative w.r.t. the radius $R$ in (\ref{eq:continuous_updateB}) and (\ref{eq:continuous_updateE}) this is not needed. 
For a function $f$ the distributional derivative in $x$-direction across a discontinuity at $x_{i-\frac{1}{2}}$ that separates two cells is given by \cite[Chapter~5]{kanwal_generalized_2004}
\begin{align}
	\partial_x f = (\partial_x f)_{\mathrm{cell}} + [f] \delta(x - x_{i-\frac{1}{2}}).
\end{align}
Here $[f] = f^+ - f^-$ denotes the jump at the discontinuity at $x_{i-\frac{1}{2}}$ and $(\partial_x f)_{\mathrm{cell}}$ 
the classical derivative within the cells left and right of the discontinuity. 
As an example, we consider the update of $E_3$ at the point value located at $(x_{i-\frac{1}{2}}, y, z)$. 
One of the parts of the integrated curl in (\ref{eq:continuous_updateE}) is given by the spherical mean of $\partial_x B_2$.
For $\partial_x B_{2}$ we have
\begin{align*}
M \{ \partial_x B^n_{2}, (x_{i-\frac{1}{2}}, y, z ); R   \} &= M \{ (\partial_x B^n_{2})_{\mathrm{cell}}, (x_{i-\frac{1}{2}}, y, z ); R   \} \\ 
&+ M \{ [B^n_{2}] \delta(x - x_{i-\frac{1}{2}}), (x_{i-\frac{1}{2}}, y, z ); R   \}.
\end{align*}
For the latter part we obtain 
\begin{equation}\label{eq:spherical_jump}
	\begin{aligned}
	&M \{ [B^n_{2}] \delta(x - x_{i-\frac{1}{2}}), (x_{i-\frac{1}{2}}, y, z ); R   \} \\ 
	&\quad = \frac{1}{4 \pi} \int_{0}^{\pi} \int_{0}^{2 \pi} [B^n_{2}]\big(x_{i-\frac{1}{2}} + R\sin(\phi)\cos(\theta), y + R\sin(\phi)\sin(\theta), \\
	&\phantom{\quad = \frac{1}{4 \pi} \int_{0}^{2 \pi} \int_{0}^{\pi} [B^n_{1}]\big((} z + R\cos(\phi)\big) \delta(R\sin(\phi)\cos(\theta)) \sin(\phi)  \,\dd\theta \,\dd\phi.
	\end{aligned}
\end{equation}
For fixed $\phi$ we get a contribution from the distribution at $\theta = \frac{\pi}{2}$ and $\theta = \frac{3\pi}{2}$. With $\sin(\phi) \neq 0$ for $\phi \in (0,\pi)$, we have 
\begin{align*}
	\delta (R\sin(\phi)\cos(\theta)) = \frac{1}{R\sin(\phi)} \left( \delta\left(\theta - \frac{\pi}{2}\right) + \delta\left(\theta - \frac{3\pi}{2}\right) \right).
\end{align*}
This gives in (\ref{eq:spherical_jump})
\begin{align}\label{eq:jump:update_x}
&\frac{1}{4 \pi R} \int_{0}^{\pi} \Big( [B^n_{2}]\big(x_{i-\frac{1}{2}}, y + R\sin(\phi), z + R\cos(\phi)\big) \nonumber \\ 
&\phantom{\frac{1}{4 \pi R} \int_{0}^{\pi} \Big( }+ [B^n_{2}]\big(x_{i-\frac{1}{2}}, y - R\sin(\phi), z + R\cos(\phi)\big) \Big) \,\dd\phi \nonumber  \\
&= \frac{1}{4 \pi R} \int_{0}^{2 \pi} \Big( [B^n_{2}]\big(x_{i-\frac{1}{2}}, y + R\sin(\phi), z + R\cos(\phi)\big) \Big) \,\dd\phi \nonumber \\
&= \frac{1}{4 \pi R} \int_{0}^{2 \pi} \Big( B^{n,+}_{2}\big(x_{i-\frac{1}{2}}, y + R\sin(\phi), z + R\cos(\phi)\big) \\
&\phantom{\frac{1}{4 \pi R} \int_{0}^{2 \pi} \Big(}-  B^{n,-}_{2}\big(x_{i-\frac{1}{2}}, y + R\sin(\phi), z + R\cos(\phi)\big) \Big) \,\dd\phi.  \nonumber
\end{align}
The remaining components are computed analogously. For completion, we give the expressions for the spherical mean of the jumps in $y$ and $z$ direction at 
discontinuities at $y_{j-\frac{1}{2}}$ and $z_{k-\frac{1}{2}}$, respectively.
We have
\begin{align}\label{eq:jump:update_y}
	&M\{ [B_1^n]\delta(y - y_{j-\frac{1}{2}}),(x, y_{j-\frac{1}{2}}, z); R \}\nonumber \\
	&\quad = \frac{1}{4 \pi R} \int_{0}^{2 \pi} \Big( B^{n,+}_{1}\big(x + R\sin(\phi), y_{j-\frac{1}{2}}, z + R\cos(\phi)\big) \\
	&\phantom{= \frac{1}{4 \pi R} \int_{0}^{2 \pi} \Big(}- B^{n,-}_{1}\big(x + R\sin(\phi), y_{j-\frac{1}{2}}, z + R\cos(\phi)\big) \Big) \,\dd\phi,\nonumber 
\end{align}
and
\begin{align}\label{eq:jump:update_z}
	&M\{ [B_2^n]\delta(z - z_{k-\frac{1}{2}}),(x, y, z_{k-\frac{1}{2}}); R \} \nonumber\\
	&\quad = \frac{1}{4 \pi R} \int_{0}^{2 \pi} \Big( B^{n,+}_{2}\big(x + R\sin(\phi), y + R\cos(\phi), z_{k-\frac{1}{2}}\big) \\
	&\phantom{= \frac{1}{4 \pi R} \int_{0}^{2 \pi} \Big(}- B^{n,-}_{2}\big(x + R\sin(\phi), y + R\cos(\phi), z_{k-\frac{1}{2}}\big) \Big) \,\dd\phi. \nonumber
\end{align}
We compute these integrals exactly similar to the approach in Section~\ref{sec:Maxwell} with SymPy beforehand by integrating cell wise reconstructions along arcs of these circles in the respective discontinuity planes.
The bound on the radius $R$ given in Section~\ref{sec:Maxwell} is still applied here. 
Thus, a necessary stability condition is CFL$\le \frac{1}{2}$.

\subsection{Initialization of the faces}
As described above, the divergence of the reconstruction reduces to the Gauss residual (\ref{eq:discrete_gauss}). 
To maintain third-order accuracy, one could initialize the facial averages with Simpson's rule.
However, this does not guarantee the divergence-free condition of this reconstruction. 
The method then carries this residual for all times, since the update (\ref{eq:stokes_time}) leaves (\ref{eq:discrete_gauss}) unchanged regardless of the value of the residual.
Thus, we initialize the facial averages from vector potentials $\bm A_0$ and $\bm C_0$
with $\bm B_0 = \nabla \times \bm A_0$ and $\bm E_0 = \nabla \times \bm C_0$ similar to what is done in constrained transport methods for MHD \cite{HELZEL20113803} and Maxwell's equations \cite{HAZRA2019298}.
The face averages of $\bm B_0$ then become the circulation of $\bm A_0$ on the boundaries of the respective faces. Analogously for $\bm E_0$ and $\bm C_0$.
The integrals are approximated with Simpson's rule. 

\subsection{The time quadrature}\label{sec:divfree_time}
Although we deviated from a classical approach to Active Flux methods for Maxwell's equations, we still want to formulate a method that uses the point values updated by the evolution operators (\ref{eq:continuous_updateB}), (\ref{eq:continuous_updateE}) and (\ref{eq:spherical_jump}) derived above.
The quadrature formula in time should satisfy three properties: First, it should be at least third-order accurate to not diminish the overall accuracy of 
the method. Second, it should not rely on boundary point values at $t_n$, since these are not uniquely defined by the applied discontinuous reconstruction. 
And third, it should use the point values at $t_{n+1}$ in order to give us diagnostics at the end of the time step.
We therefore use the two-point right-sided Gauss--Radau rule, referred to as Radau 2 in the following.
We have
\begin{align}\label{eq:radau}
	\frac{1}{\Delta t} \int_{t_n}^{t_{n+1}} f(t) \,\dd t \approx \frac{3}{4} f\!\left(t_n + \frac{1}{3}\Delta t\right)
	+ \frac{1}{4} f(t_{n+1}),
\end{align}
with nodes $\tau_1 = \frac{1}{3}$, $\tau_2 = 1$ and weights $w_1 = \frac{3}{4}$, $w_2 = \frac{1}{4}$. 

One time step of the method is summarized in Algorithm~\ref{alg:divfree}.
\begin{algorithm}[htbp]
\caption{Divergence-free method.}
\label{alg:divfree}
\begin{algorithmic}[1]
\State \textbf{Reconstruct} $\bm B^{rec}$, $\bm E^{rec}$ with (\ref{eq:divfree_ansatz}) on every cell.
\State \textbf{Evolve} the nine point values of every face to the Radau times with $R = c \tau\Delta t$ using the evolution operator (\ref{eq:continuous_updateB}) and (\ref{eq:continuous_updateE}) with distributional curl (\ref{eq:spherical_jump}).
\State \textbf{Edge integrals:} Compute the edge integrals $C^1_{j\pm\frac{1}{2}, k\pm\frac{1}{2}}$, $C^2_{i\pm\frac{1}{2}, k\pm\frac{1}{2}}$, $C^3_{i\pm\frac{1}{2}, j\pm\frac{1}{2}}$ 
with Simpson's rule, for every step of the Radau rule (\ref{eq:radau}) in time.
\State \textbf{Update facial averages by Stokes} (\ref{eq:stokes_time}), (\ref{eq:stokes_time_yz}) on the respective faces.
\State \textbf{Update the five remaining facial moments} for $\bm B$ and $\bm E$ by (\ref{eq:facial_FD}) from the evolved face point values to build the reconstruction in the following time step.
\end{algorithmic}
\end{algorithm}

\section{Numerical examples}\label{sec:numerics}

In this section we consider numerical examples for the proposed methods. For brevity we refer to the
method of Section~\ref{sec:Maxwell} as the \textit{classical Active Flux method} and to the method of
Section~\ref{sec:divfree} as the \textit{divergence-free method}. 
In all examples we set $c=1$. 

We measure the errors $e_{i,j,k}$ for the method of Section~\ref{sec:Maxwell} by comparing the cell averages with cell averages computed with Simpson's rule from the exact solution, and the method of Section~\ref{sec:divfree} by
comparing the face-averages with face-averages computed with Simpson's rule from the exact solution.
We consider the norms 
\begin{align*}
	\Vert e \Vert_{L^1} &= \sum_{i,j,k} \vert e_{i,j,k}\vert \, \Delta x \Delta y \Delta z, \qquad
	\Vert e \Vert_{L^2} = \Big( \sum_{i,j,k} e_{i,j,k}^2 \, \Delta x \Delta y \Delta z \Big)^{1/2}, \\
	\Vert e \Vert_{L^\infty} &= \max_{i,j,k} \vert e_{i,j,k}\vert.
\end{align*}
For the face averages, $e_{i,j,k}$ denotes the error associated with the respective component on the face, e.g. on the $x$-faces for $B_1$ and $E_1$. The sums run over all of these faces.
To compare both methods, we compare the divergence and the total energy. 
For the method of Section~\ref{sec:Maxwell}, one measures the divergence of the reconstruction in the $L^2$
norm, in this case the result is given by
\begin{align}
    \Vert\nabla \cdot \bm{B}(\cdot, t_n)\Vert_{L^2}^2
    \approx \sum_{i,j,k} \int_{\mathcal{C}_{i,j,k}} \left(\partial_x B_{1}^{rec,n}
    + \partial_y B_{2}^{rec,n}
    + \partial_z B_{3}^{rec,n}\right)^2 \,\dd x \,\dd y \,\dd z. \label{eq:Rec_divergence}
\end{align}
Since the divergence of the reconstruction of the method described in Section~\ref{sec:divfree_rec} is constant within a cell, the respective quantity 
reduces to 
\begin{align}
	\Vert\nabla \cdot \bm{B}(\cdot, t_n)\Vert_{L^2}^2 \approx \sum_{i,j,k} \Delta x \Delta y \Delta z \left(\partial_x B_{1}^{rec,n}
    + \partial_y B_{2}^{rec,n}
    + \partial_z B_{3}^{rec,n}\right)^2.
\end{align}
We also monitor the maximum value of the divergence of each reconstruction. For the divergence-free method 
from Section~\ref{sec:divfree_rec} this quantity is constant in each cell and therefore the maximum value is straightforward to 
compute. 
For the classical Active Flux method from Section~\ref{sec:Maxwell}, we evaluate the divergence of the cell's 
reconstruction at the 26 points on the boundary and its center and take the maximum.
Furthermore, we present the numerical results of the total electromagnetic energy given by
\begin{align}
	\mathcal{E}(t) = \frac{1}{2} \int \left( |\bm{B}|^2 + |\bm{E}|^2\right)\,\dd\bm x.
\end{align} 
For the method of Section~\ref{sec:Maxwell} we compute the energy directly from the cell averages, and for the method of Section~\ref{sec:divfree}
 we compute the energy from the average of the reconstruction over the cell.
We monitor the relative energy error
\begin{align}\label{eq:energy_error}
    \frac{\mathcal{E}(t) - \mathcal{E}(0)}{\mathcal{E}(0)}.
\end{align}

\subsection{Plane wave}\label{sec:planewaveBalsara}
The first example is a plane wave from \cite{BALSARA20191108} on the domain $[-\frac{1}{2},\frac{1}{2}]^3$ with periodic boundary conditions. 
There the authors considered a plane polarized wave. The magnetic field is given by the vector potential 
\begin{equation}
	\bm{A} = \frac{1}{2 \pi} \sin\left( 2\pi(x + y - \sqrt{2}t)\right) \begin{pmatrix} 0 \\ 1 \\ 0 \end{pmatrix},
\end{equation}
and the electric field is given by the vector potential
\begin{equation}
	\bm{C} = -\frac{1}{2 \sqrt{2} \pi} \sin \left(2\pi(x + y - \sqrt{2}t)\right) \begin{pmatrix} 0 \\ 0 \\ 1 \end{pmatrix}.
\end{equation}
These give the magnetic and electric fields
\begin{align}
	\bm{B} &= \nabla \times \bm{A} = \cos\left( 2\pi(x + y - \sqrt{2}t)\right) \begin{pmatrix} 0 \\ 0 \\ 1 \end{pmatrix}, \qquad \nonumber \\
	\bm{E} &= \nabla \times \bm{C} = \frac{1}{\sqrt{2}} \cos\left( 2\pi(x + y - \sqrt{2}t)\right) \begin{pmatrix} -1 \\ 1 \\ 0 \end{pmatrix}.
\end{align}
For this example, the Simpson initialization of the face average gives a vanishing Gauss residual.
In Tables~\ref{tab:eoc_planewave_classical} and~\ref{tab:eoc_planewave_divfree} we present the errors of both methods over one period, 
i.e., until $t = \frac{1}{\sqrt{2}} \approx 0.707$.

\begin{table}[htbp]
\centering
\small
\setlength{\tabcolsep}{4pt}
\adjustbox{max width=\linewidth}{%
\begin{tabular}{llrcrcrc}
\toprule
\textbf{Comp.} & \textbf{$N^3$} & $L^1$\textbf{-Error} & \textbf{EOC} & $L^2$\textbf{-Error} & \textbf{EOC} & $L^\infty$\textbf{-Error} & \textbf{EOC} \\
\midrule
  $B_3$ & $16^3$ & $1.272 \cdot 10^{-3}$ & --- & $1.431 \cdot 10^{-3}$ & --- & $2.024 \cdot 10^{-3}$ & --- \\
   & $32^3$ & $1.667 \cdot 10^{-4}$ & 2.93 & $1.857 \cdot 10^{-4}$ & 2.95 & $2.627 \cdot 10^{-4}$ & 2.95 \\
   & $64^3$ & $2.113 \cdot 10^{-5}$ & 2.98 & $2.348 \cdot 10^{-5}$ & 2.98 & $3.321 \cdot 10^{-5}$ & 2.98 \\
   & $128^3$ & $2.645 \cdot 10^{-6}$ & 3.00 & $2.939 \cdot 10^{-6}$ & 3.00 & $4.156 \cdot 10^{-6}$ & 3.00 \\
\midrule
  $E_2$ & $16^3$ & $9.014 \cdot 10^{-4}$ & --- & $1.014 \cdot 10^{-3}$ & --- & $1.434 \cdot 10^{-3}$ & --- \\
   & $32^3$ & $1.180 \cdot 10^{-4}$ & 2.93 & $1.315 \cdot 10^{-4}$ & 2.95 & $1.860 \cdot 10^{-4}$ & 2.95 \\
   & $64^3$ & $1.494 \cdot 10^{-5}$ & 2.98 & $1.660 \cdot 10^{-5}$ & 2.99 & $2.348 \cdot 10^{-5}$ & 2.99 \\
   & $128^3$ & $1.871 \cdot 10^{-6}$ & 3.00 & $2.079 \cdot 10^{-6}$ & 3.00 & $2.940 \cdot 10^{-6}$ & 3.00 \\
\bottomrule
\end{tabular}}
\caption{Errors and experimental orders of convergence for the classical Active Flux method of
Section~\ref{sec:Maxwell} on the plane wave of Section~\ref{sec:planewaveBalsara} at $t = \frac{1}{\sqrt{2}}$ and $\text{CFL} = 0.45$.}
\label{tab:eoc_planewave_classical}
\end{table}

\begin{table}[htbp]
\centering
\small
\setlength{\tabcolsep}{4pt}
\adjustbox{max width=\linewidth}{%
\begin{tabular}{llrcrcrc}
\toprule
\textbf{Comp.} & \textbf{$N^3$} & $L^1$\textbf{-Error} & \textbf{EOC} & $L^2$\textbf{-Error} & \textbf{EOC} & $L^\infty$\textbf{-Error} & \textbf{EOC} \\
\midrule
  $B_3$ & $16^3$ & $2.645 \cdot 10^{-3}$ & --- & $2.971 \cdot 10^{-3}$ & --- & $4.202 \cdot 10^{-3}$ & --- \\
   & $32^3$ & $3.480 \cdot 10^{-4}$ & 2.93 & $3.876 \cdot 10^{-4}$ & 2.94 & $5.482 \cdot 10^{-4}$ & 2.94 \\
   & $64^3$ & $4.420 \cdot 10^{-5}$ & 2.98 & $4.913 \cdot 10^{-5}$ & 2.98 & $6.948 \cdot 10^{-5}$ & 2.98 \\
   & $128^3$ & $5.549 \cdot 10^{-6}$ & 2.99 & $6.165 \cdot 10^{-6}$ & 2.99 & $8.718 \cdot 10^{-6}$ & 2.99 \\
\midrule
  $E_2$ & $16^3$ & $1.937 \cdot 10^{-3}$ & --- & $2.138 \cdot 10^{-3}$ & --- & $2.970 \cdot 10^{-3}$ & --- \\
   & $32^3$ & $2.494 \cdot 10^{-4}$ & 2.96 & $2.766 \cdot 10^{-4}$ & 2.95 & $3.895 \cdot 10^{-4}$ & 2.93 \\
   & $64^3$ & $3.145 \cdot 10^{-5}$ & 2.99 & $3.492 \cdot 10^{-5}$ & 2.99 & $4.933 \cdot 10^{-5}$ & 2.98 \\
   & $128^3$ & $3.929 \cdot 10^{-6}$ & 3.00 & $4.363 \cdot 10^{-6}$ & 3.00 & $6.169 \cdot 10^{-6}$ & 3.00 \\
\bottomrule
\end{tabular}}
\caption{Errors of the face averages and experimental orders of convergence for the divergence-free method of
Section~\ref{sec:divfree} on the plane wave of Section~\ref{sec:planewaveBalsara} at $t = \frac{1}{\sqrt{2}}$ and $\text{CFL} = 0.45$.}
\label{tab:eoc_planewave_divfree}
\end{table}
One can see that both methods are of third order in all considered norms. For the magnetic field we see that both methods perform better than the third-order DGTD method of \cite{BALSARA20191108}. 
We do not compare the errors of the electric field, since \cite{BALSARA20191108} reports the electric displacement $\bm D$ in SI units, with errors computed at the cell faces. 
On all considered grids, the error observed for the method from Section~\ref{sec:divfree} is roughly twice as large as the error observed for the method from Section~\ref{sec:Maxwell}.

\subsection{Stability and long-time behavior}\label{sec:stability}
\begin{figure}[htbp]
	\centering
	\includegraphics[width=0.49\textwidth]{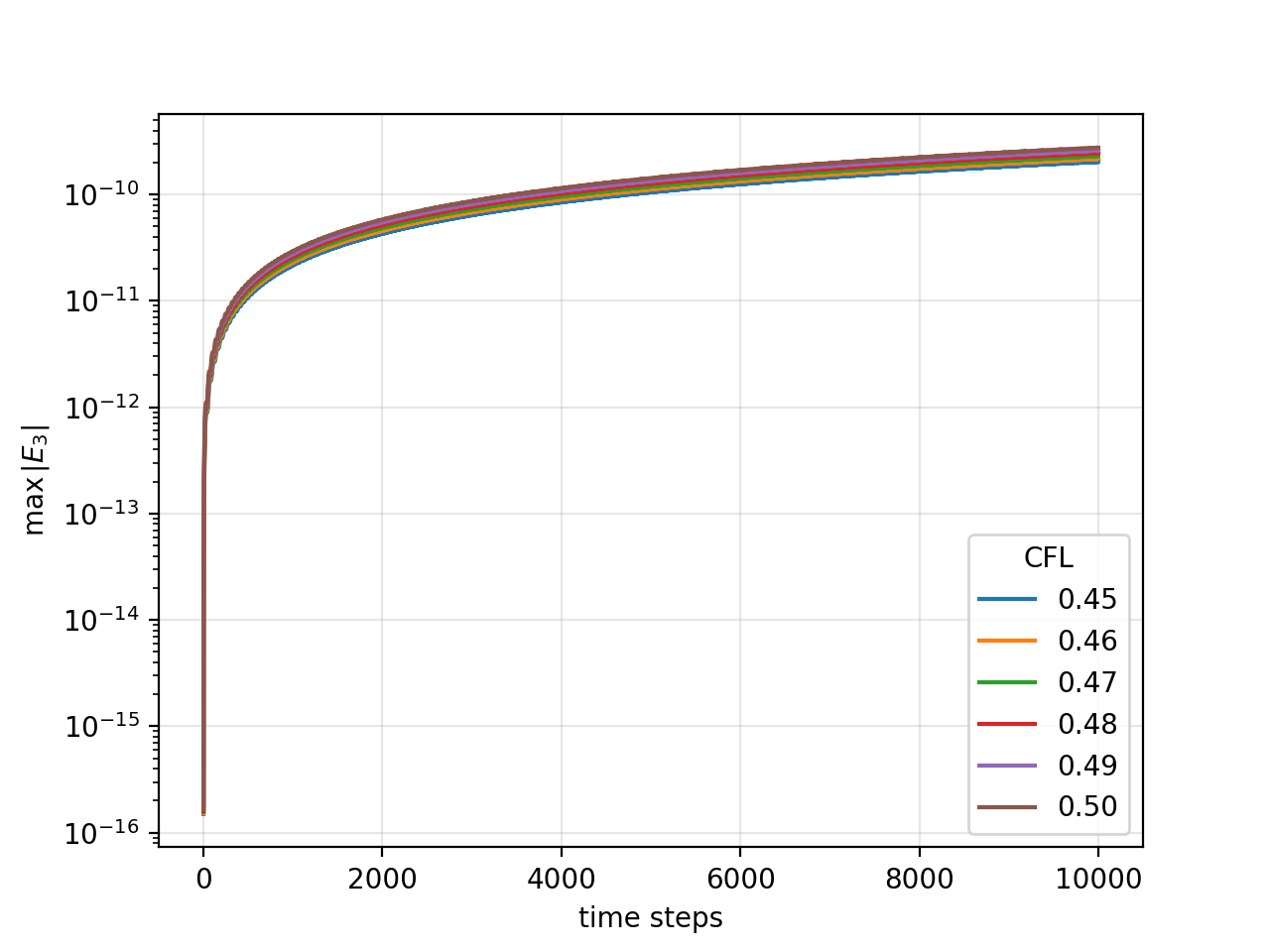}
	\includegraphics[width=0.49\textwidth]{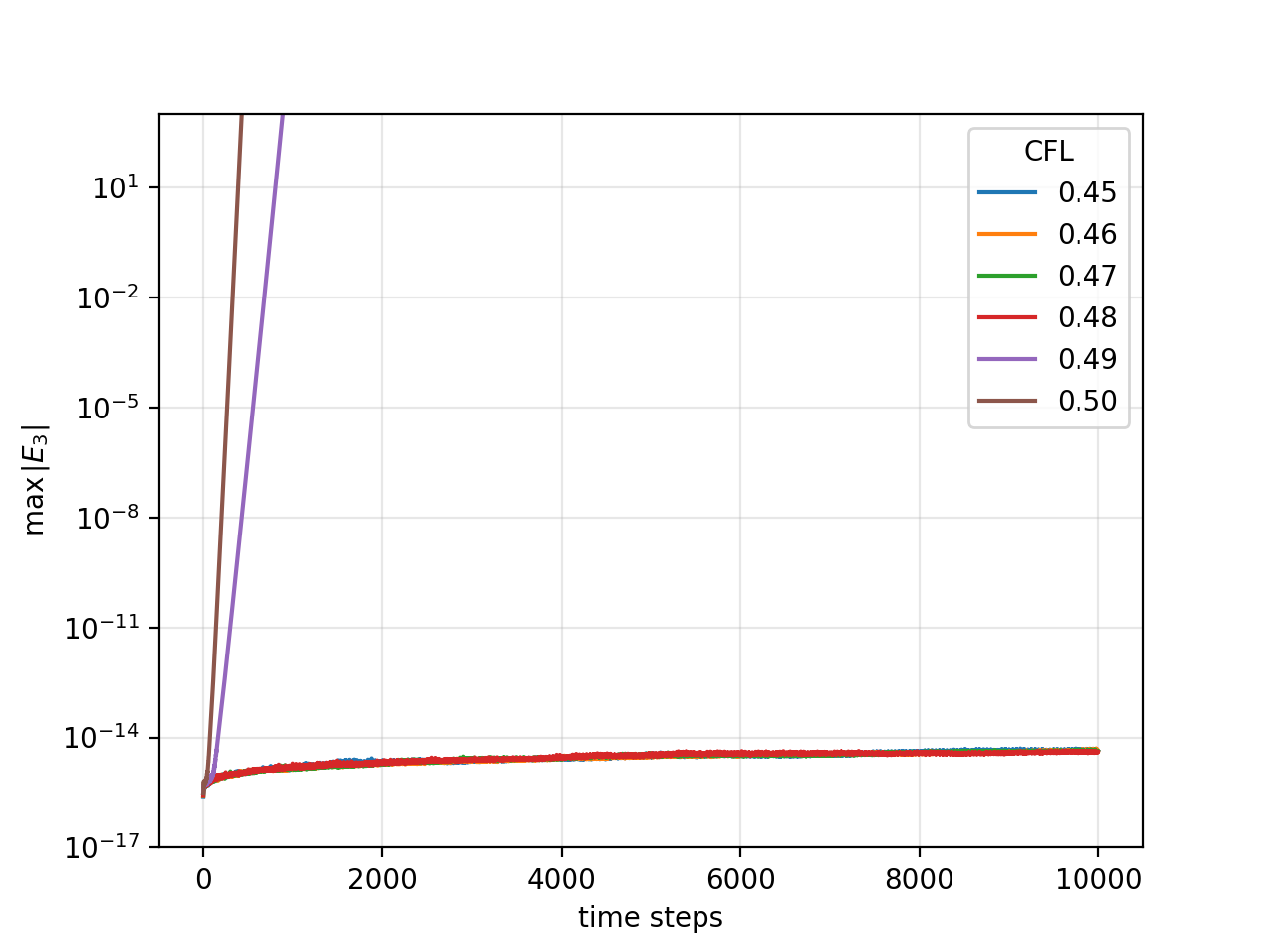 }
	\caption{Maximum of $\vert E_3 \vert $ of the cell averages for the classical (left) and the face averages for the divergence-free method (right) over $10000$ time steps 
	for the example (\ref{eq:plane_wave2}) on a $32^3$ grid, for the divergence-free method for CFL numbers from $0.45$ to $0.5$.}
	\label{fig:cfl_edges}
\end{figure}

Next, we want to investigate the stability and long-time behavior of both methods. 
For this, we consider an extension of the previous example. The solution is given by
\begin{align}\label{eq:plane_wave2}
  \bm E(\bm x,t) &= \frac{1}{\sqrt2} \cos\left(\pi(x+y+z-\sqrt3 t)\right) \begin{pmatrix}1\\-1\\0\end{pmatrix}, \\
  \bm B(\bm x,t) &= \frac{1}{\sqrt6} \cos\left(\pi(x+y+z-\sqrt3 t)\right) \begin{pmatrix}1\\1\\-2\end{pmatrix}. \nonumber
\end{align}
Here we do not need an initialization of the facial averages from the potential, because for this wave Simpson's rule is enough to guarantee that 
the discrete Gauss' law is satisfied. We consider the domain $[-1,1]^3$ and periodic boundary conditions.
We are first interested in investigating the stability. To monitor the behavior of the methods, we consider the approximation of 
$E_3$. A stable method should keep the maximum of that component small. 
We found that the presented divergence-free method does not reach the same stability as the classical method.
In Figure~\ref{fig:cfl_edges} we present the maximum of $\vert E_3 \vert $ of the cell averages for the classical method 
and the face averages for the divergence-free method over $10000$ time steps for CFL numbers from $0.45$ to $0.5$. 
One can see that the classical method is stable for all CFL numbers considered. Its error grows only slowly and stays at about $10^{-10}$.
The divergence-free method appears to be stable up to a CFL number of $0.48$ in this example, while the error grows exponentially for CFL $0.49$ and $0.5$.

\begin{figure}[htbp]
	\centering
	\includegraphics[width=\textwidth]{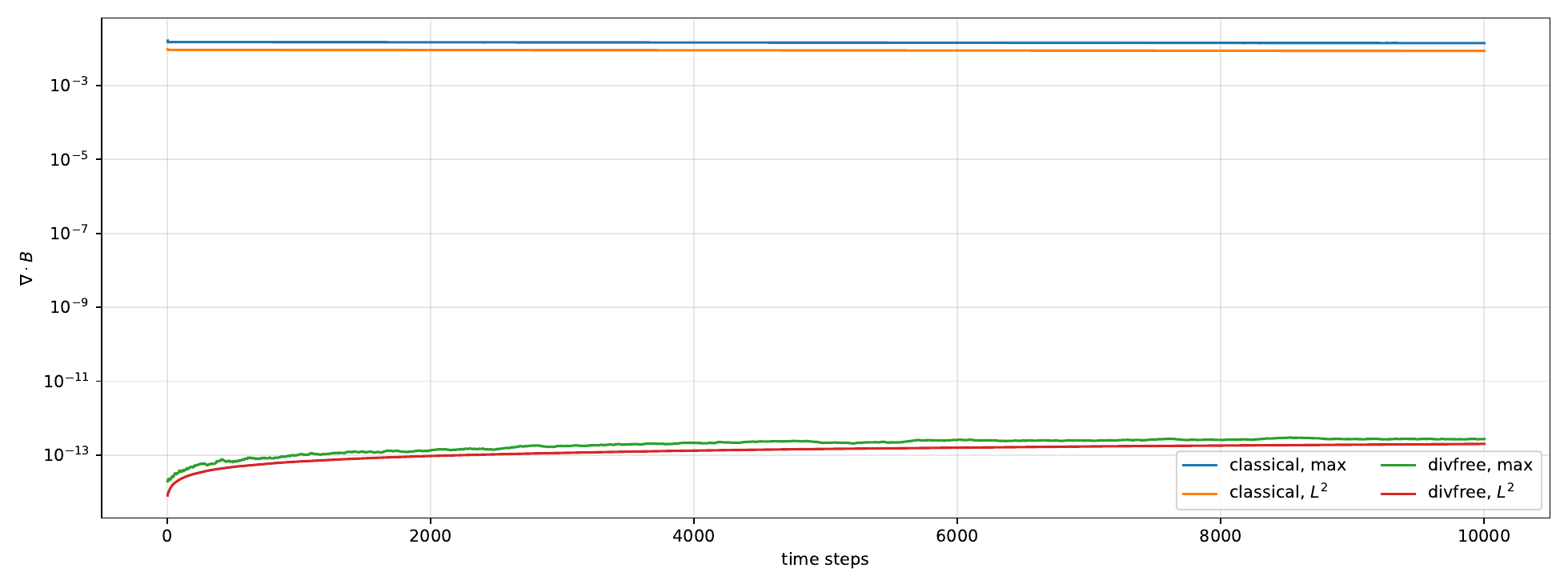}
	\caption{Maximum and $L^2$ norm of $\nabla \cdot \bm B$ over $10000$ time steps, $N^3 = 32^3$, with $\text{CFL} = 0.45$ for both methods.}
	\label{fig:longtime_div}
\end{figure}

\begin{figure}[htbp]
	\centering
	\includegraphics[width=0.49\textwidth]{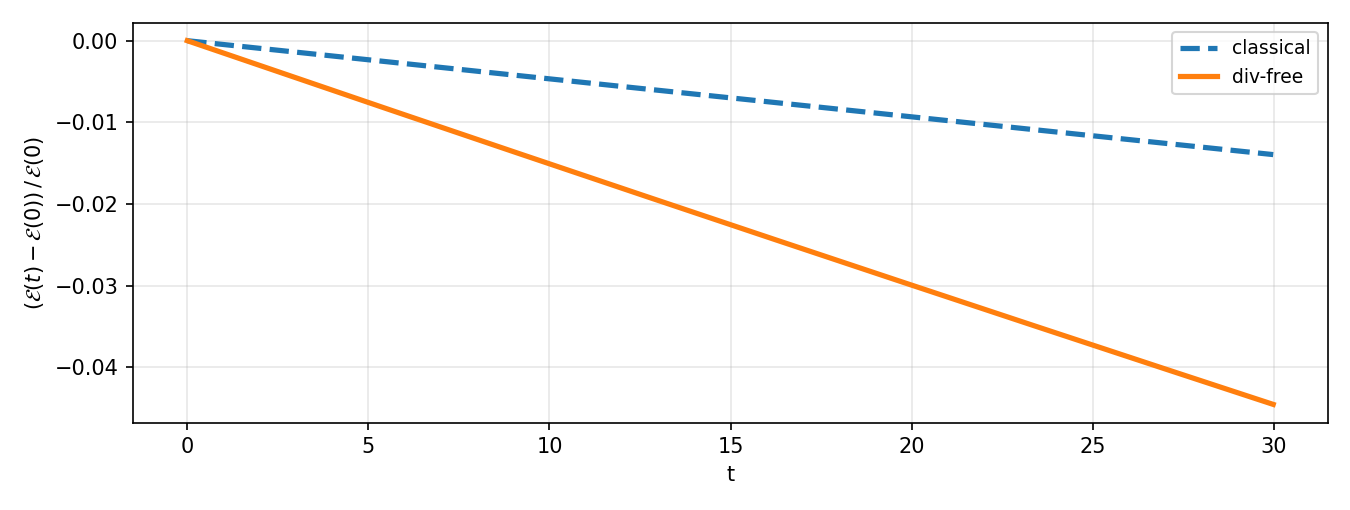}
	\includegraphics[width=0.49\textwidth]{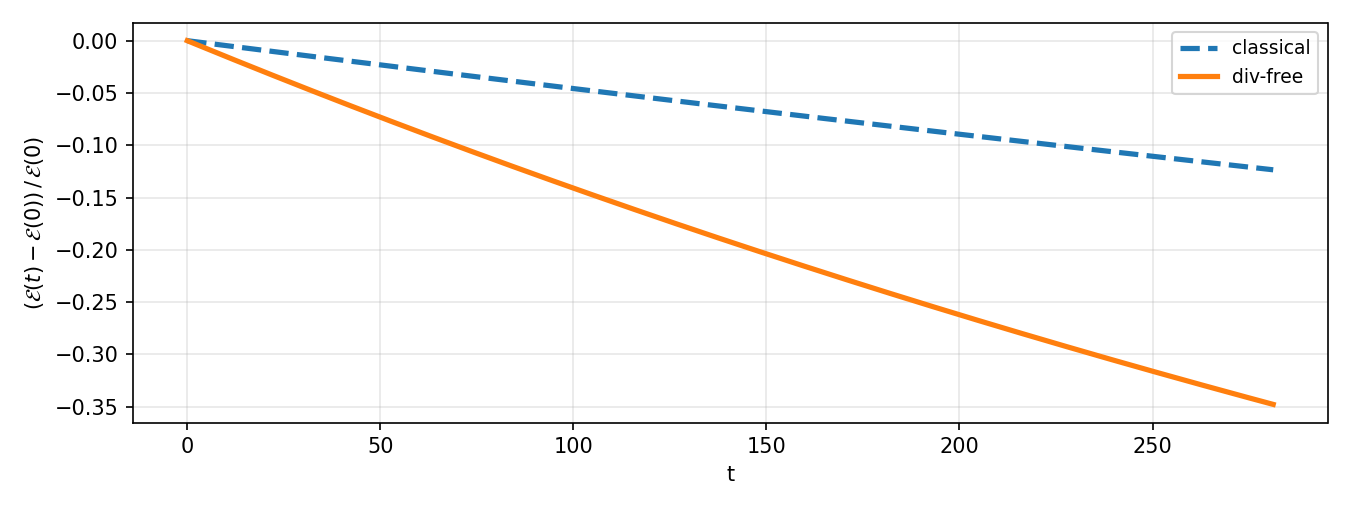}
	\caption{Relative energy error for grid with $N^3 = 32^3$ cells and $\text{CFL} = 0.45$ for both methods. Left: $t \in [0,30]$, right: over all $10000$ time steps, i.e.\ $t \in [0,281.25]$.}
	\label{fig:longtime_energy}
\end{figure}

We also considered the long-time behavior of the proposed methods and investigated the preservation of the divergence and the total energy for the initial data above. 
In Figure~\ref{fig:longtime_div} we present the maximum and the $L^2$ norm of the divergence of both methods and in 
Figure~\ref{fig:longtime_energy} we present the relative energy (\ref{eq:energy_error}) for both methods.
The divergence of $\bm E$ behaves in the same way and is not shown. The classical Active Flux method keeps its divergence at the level
set by the initial data. The divergence-free method holds its
divergence at round-off over the whole run. However, it loses more energy than the classical method. 

\subsection{Oblique wave}\label{sec:oblique}
\begin{figure}[htbp]
	\centering
	\begin{subfigure}{\textwidth}
		\includegraphics[width=0.49\linewidth]{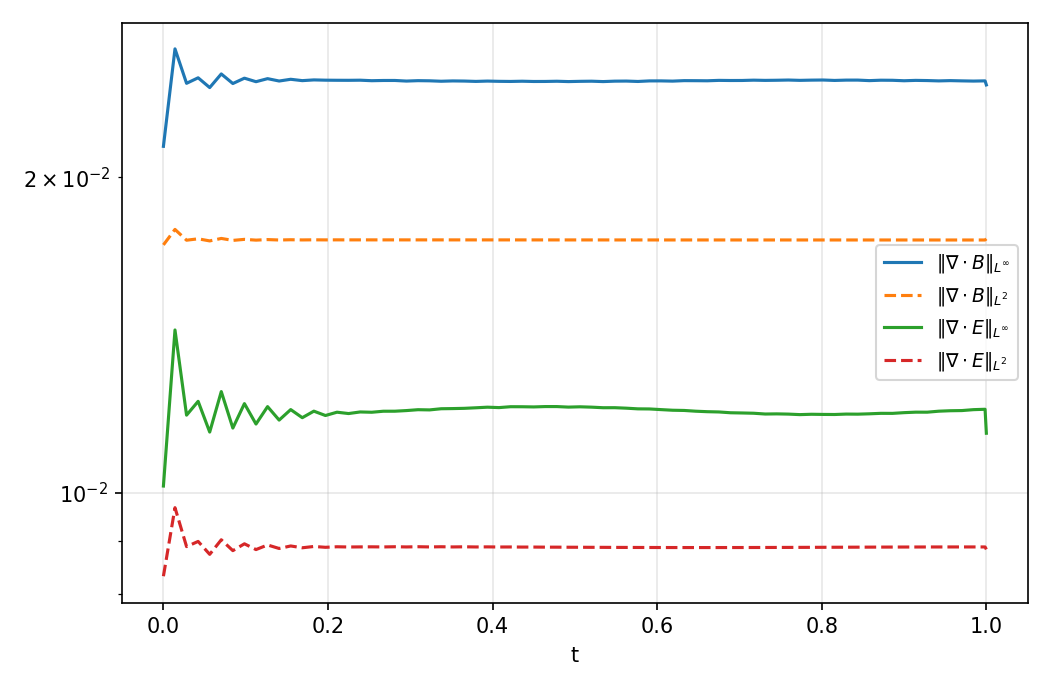}
		\includegraphics[width=0.49\linewidth]{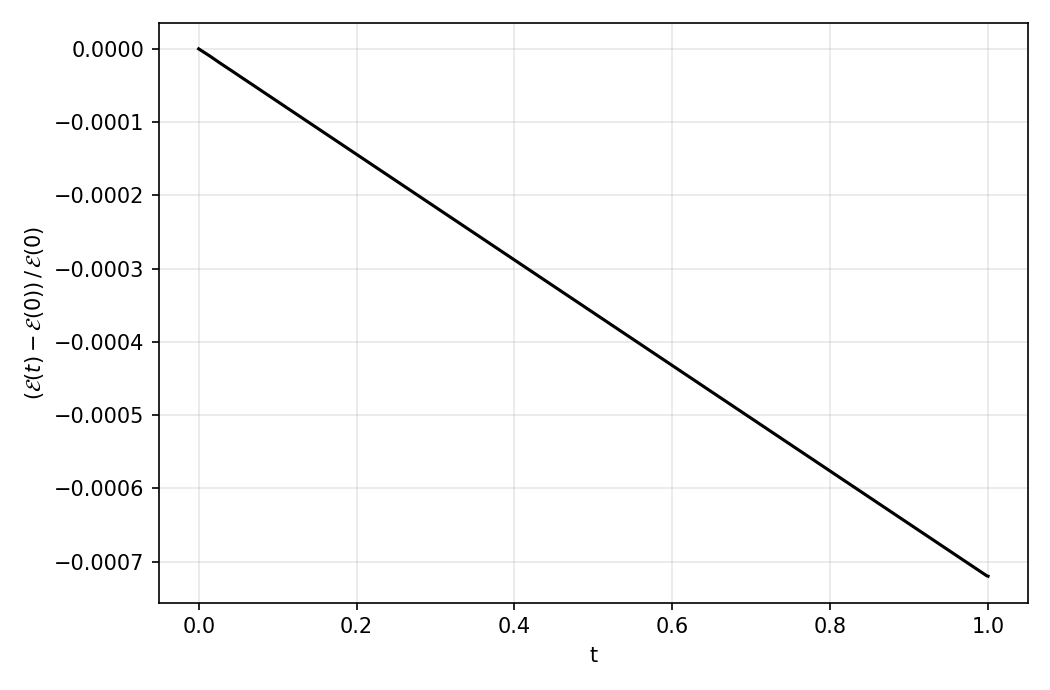}
		\caption{Classical Active Flux method of Section~\ref{sec:Maxwell}.}
	\end{subfigure}

	\begin{subfigure}{\textwidth}
		\includegraphics[width=0.49\linewidth]{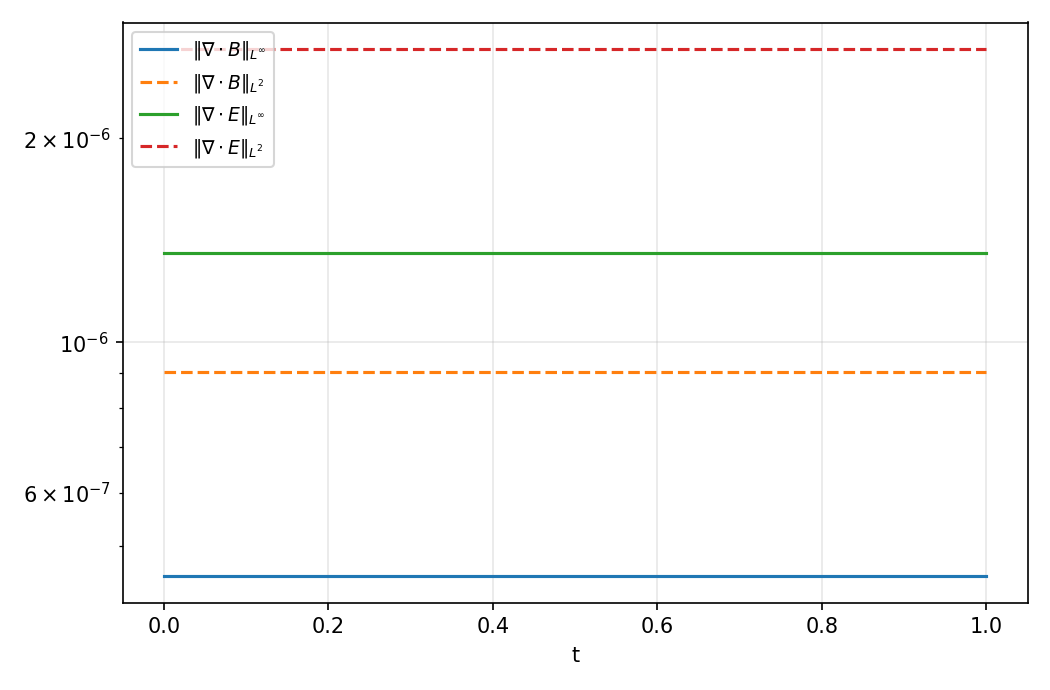}
		\includegraphics[width=0.49\linewidth]{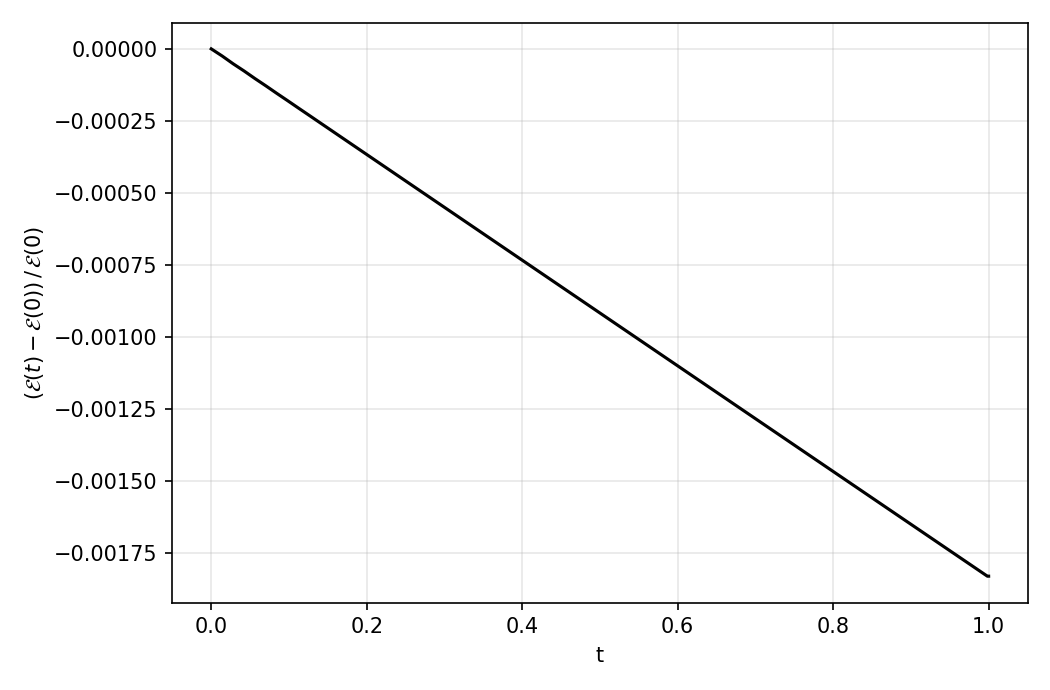}
		\caption{Divergence-free method of Section~\ref{sec:divfree} with face averages initialized by Simpson's rule.}
	\end{subfigure}

	\begin{subfigure}{\textwidth}
		\includegraphics[width=0.49\linewidth]{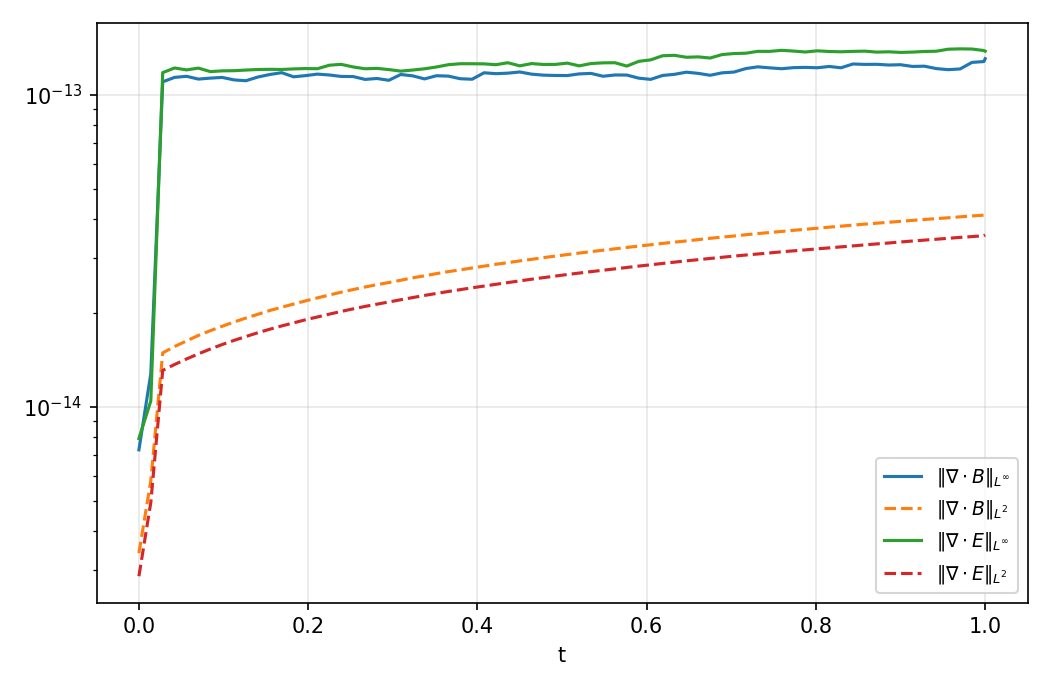}
		\includegraphics[width=0.49\linewidth]{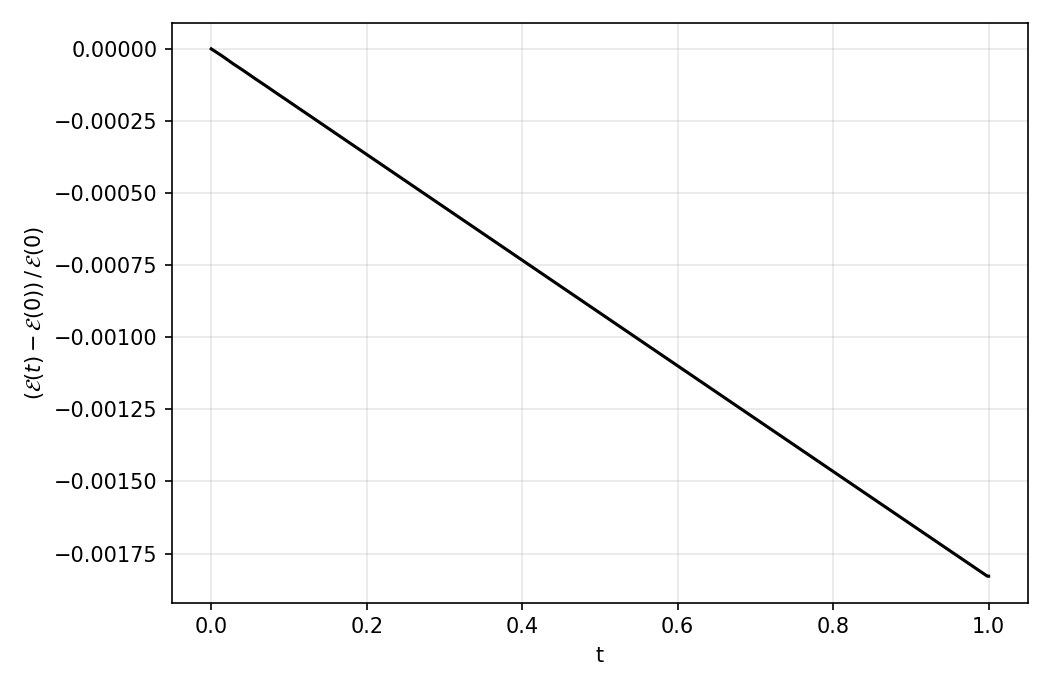}
		\caption{Divergence-free method of Section~\ref{sec:divfree} with face averages initialized from the potential.}
	\end{subfigure}
	\caption{Divergence of $\bm B$ and $\bm E$ in the $L^\infty$ and $L^2$ norm (left) and relative energy error (\ref{eq:energy_error}) (right) over $t \in [0,1]$ 
	for the oblique plane wave (\ref{eq:oblique_wave}) with $N^3 = 64^3$ cells and $\text{CFL} = 0.45$ for all methods.}
	\label{fig:planewave_divenergy}
\end{figure}

Another example to investigate both methods is the following oblique wave with wave vector $\bm k = \pi(1,2,2)^T$
\begin{align}\label{eq:oblique_wave}
  \bm E(\bm x,t) &= \frac{1}{\sqrt5} \cos\left(\pi(x+2y+2z-3t)\right) \begin{pmatrix}2\\-1\\0\end{pmatrix}, \\
  \bm B(\bm x,t) &= \frac{1}{3\sqrt5} \cos\left(\pi(x+2y+2z-3t)\right) \begin{pmatrix}2\\4\\-5\end{pmatrix}. \nonumber
\end{align}
We again consider the domain $[-1,1]^3$ with periodic boundary conditions. For the initialization 
of the facial averages we use the following vector potentials
\begin{align}
  \bm A_0(\bm x) &= \frac{1}{3\sqrt5\pi} \sin\left(\pi(x+2y+2z)\right) \begin{pmatrix}2\\-1\\0\end{pmatrix}, \\
  \bm C_0(\bm x) &= -\frac{1}{9\sqrt5\pi} \sin\left(\pi(x+2y+2z)\right) \begin{pmatrix}2\\4\\-5\end{pmatrix}, \nonumber
\end{align}
so that $\nabla\times\bm A_0 = \bm B(\bm x, 0)$ and $\nabla\times\bm C_0 = \bm E(\bm x, 0)$. 
In the tables below we present the order of convergence for both methods and for their respective divergences as well as the 
relative energy error (\ref{eq:energy_error}) at $t=1$. 
For all computations we consider a CFL number of $0.45$ for the classical Active Flux method and the divergence-free method.
In Tables~\ref{tab:eoc_oblique_classical} and~\ref{tab:eoc_oblique_divfree} we present the accuracy of the classical Active Flux method and the divergence-free method for the oblique wave (\ref{eq:oblique_wave}) at $t = 1$.
Both methods confirm third-order accuracy in all components. 
The observed errors are again about twice as large for the method from Section~\ref{sec:divfree}. An exception is $E_3$, which vanishes in the exact solution. 
For this component the error of the divergence-free method is smaller.
In Tables~\ref{tab:eoc_div} and~\ref{tab:eoc_energy} we distinguish between the classical method, the divergence-free method with Simpson's 
initialization of the facial averages and the divergence-free method with initialization from the vector potentials.
We see that the divergence of the classical method converges with second order and the divergence-free method with Simpson's initialization 
converges with fourth order. The reason for these discrepancies is that the divergence of the classical method is computed 
from the derivative of the reconstruction and thus loses one order of accuracy. The divergence-free method shows this behavior
because the faces are initialized with Simpson's rule, which is a fourth-order method. This order gets carried over throughout 
the method. 
In Figure~\ref{fig:planewave_divenergy} we present plots of the divergence and the relative energy error for three configurations. 
It can be seen that the divergence remains bounded for all methods. For the classical method, it settles slightly above its initial level after a few time steps. 
For the divergence-free method with Simpson initialisation, it remains at the level set by the initialisation. 
With initialisation from the vector potentials, it rises during the first few time steps to approximately $10^{-13}$ in the maximum norm and thus remains at round-off.
For this example the divergence-free method needs the initialization from the vector potentials to keep the divergence at round-off. 
The energy dissipates for all methods, but the classical Active Flux method dissipates less energy than the divergence-free method. The dissipation 
of the divergence-free method is independent of the initialization of the facial averages.

\begin{table}[htbp]
\centering
\small
\setlength{\tabcolsep}{4pt}
\adjustbox{max width=\linewidth}{%
\begin{tabular}{llrcrcrc}
\toprule
\textbf{Comp.} & \textbf{$N^3$} & $L^1$\textbf{-Error} & \textbf{EOC} & $L^2$\textbf{-Error} & \textbf{EOC} & $L^\infty$\textbf{-Error} & \textbf{EOC} \\
\midrule
  $B_1$ & $16^3$ & $2.988 \cdot 10^{-2}$ & --- & $1.166 \cdot 10^{-2}$ & --- & $5.719 \cdot 10^{-3}$ & --- \\
   & $32^3$ & $4.223 \cdot 10^{-3}$ & 2.82 & $1.656 \cdot 10^{-3}$ & 2.82 & $8.260 \cdot 10^{-4}$ & 2.79 \\
   & $64^3$ & $5.449 \cdot 10^{-4}$ & 2.95 & $2.140 \cdot 10^{-4}$ & 2.95 & $1.070 \cdot 10^{-4}$ & 2.95 \\
   & $128^3$ & $6.876 \cdot 10^{-5}$ & 2.99 & $2.701 \cdot 10^{-5}$ & 2.99 & $1.350 \cdot 10^{-5}$ & 2.99 \\
\midrule
  $B_2$ & $16^3$ & $5.915 \cdot 10^{-2}$ & --- & $2.308 \cdot 10^{-2}$ & --- & $1.133 \cdot 10^{-2}$ & --- \\
   & $32^3$ & $8.435 \cdot 10^{-3}$ & 2.81 & $3.308 \cdot 10^{-3}$ & 2.80 & $1.650 \cdot 10^{-3}$ & 2.78 \\
   & $64^3$ & $1.091 \cdot 10^{-3}$ & 2.95 & $4.283 \cdot 10^{-4}$ & 2.95 & $2.141 \cdot 10^{-4}$ & 2.95 \\
   & $128^3$ & $1.376 \cdot 10^{-4}$ & 2.99 & $5.405 \cdot 10^{-5}$ & 2.99 & $2.702 \cdot 10^{-5}$ & 2.99 \\
\midrule
  $B_3$ & $16^3$ & $7.367 \cdot 10^{-2}$ & --- & $2.876 \cdot 10^{-2}$ & --- & $1.418 \cdot 10^{-2}$ & --- \\
   & $32^3$ & $1.053 \cdot 10^{-2}$ & 2.81 & $4.127 \cdot 10^{-3}$ & 2.80 & $2.056 \cdot 10^{-3}$ & 2.79 \\
   & $64^3$ & $1.362 \cdot 10^{-3}$ & 2.95 & $5.348 \cdot 10^{-4}$ & 2.95 & $2.672 \cdot 10^{-4}$ & 2.94 \\
   & $128^3$ & $1.719 \cdot 10^{-4}$ & 2.99 & $6.752 \cdot 10^{-5}$ & 2.99 & $3.376 \cdot 10^{-5}$ & 2.98 \\
\midrule
  $E_1$ & $16^3$ & $8.963 \cdot 10^{-2}$ & --- & $3.497 \cdot 10^{-2}$ & --- & $1.716 \cdot 10^{-2}$ & --- \\
   & $32^3$ & $1.267 \cdot 10^{-2}$ & 2.82 & $4.969 \cdot 10^{-3}$ & 2.82 & $2.478 \cdot 10^{-3}$ & 2.79 \\
   & $64^3$ & $1.635 \cdot 10^{-3}$ & 2.95 & $6.419 \cdot 10^{-4}$ & 2.95 & $3.209 \cdot 10^{-4}$ & 2.95 \\
   & $128^3$ & $2.063 \cdot 10^{-4}$ & 2.99 & $8.102 \cdot 10^{-5}$ & 2.99 & $4.051 \cdot 10^{-5}$ & 2.99 \\
\midrule
  $E_2$ & $16^3$ & $4.391 \cdot 10^{-2}$ & --- & $1.719 \cdot 10^{-2}$ & --- & $8.542 \cdot 10^{-3}$ & --- \\
   & $32^3$ & $6.297 \cdot 10^{-3}$ & 2.80 & $2.477 \cdot 10^{-3}$ & 2.80 & $1.238 \cdot 10^{-3}$ & 2.79 \\
   & $64^3$ & $8.189 \cdot 10^{-4}$ & 2.94 & $3.215 \cdot 10^{-4}$ & 2.95 & $1.606 \cdot 10^{-4}$ & 2.95 \\
   & $128^3$ & $1.035 \cdot 10^{-4}$ & 2.98 & $4.063 \cdot 10^{-5}$ & 2.98 & $2.031 \cdot 10^{-5}$ & 2.98 \\
\midrule
  $E_3$ & $16^3$ & $3.045 \cdot 10^{-3}$ & --- & $1.197 \cdot 10^{-3}$ & --- & $5.965 \cdot 10^{-4}$ & --- \\
   & $32^3$ & $5.384 \cdot 10^{-4}$ & 2.50 & $2.113 \cdot 10^{-4}$ & 2.50 & $1.055 \cdot 10^{-4}$ & 2.50 \\
   & $64^3$ & $7.235 \cdot 10^{-5}$ & 2.90 & $2.841 \cdot 10^{-5}$ & 2.89 & $1.420 \cdot 10^{-5}$ & 2.89 \\
   & $128^3$ & $9.187 \cdot 10^{-6}$ & 2.98 & $3.607 \cdot 10^{-6}$ & 2.98 & $1.804 \cdot 10^{-6}$ & 2.98 \\
\bottomrule
\end{tabular}}
\caption{Errors and experimental orders of convergence for the method of Section~\ref{sec:Maxwell} on the
oblique plane wave (\ref{eq:oblique_wave}) at $t=1$ and $\text{CFL} = 0.45$.}
\label{tab:eoc_oblique_classical}
\end{table}

\begin{table}[htbp]
\centering
\small
\setlength{\tabcolsep}{4pt}
\adjustbox{max width=\linewidth}{%
\begin{tabular}{llrcrcrc}
\toprule
\textbf{Comp.} & \textbf{$N^3$} & $L^1$\textbf{-Error} & \textbf{EOC} & $L^2$\textbf{-Error} & \textbf{EOC} & $L^\infty$\textbf{-Error} & \textbf{EOC} \\
\midrule
  $B_1$ & $16^3$ & $7.279 \cdot 10^{-2}$ & --- & $2.877 \cdot 10^{-2}$ & --- & $1.438 \cdot 10^{-2}$ & --- \\
   & $32^3$ & $1.054 \cdot 10^{-2}$ & 2.79 & $4.151 \cdot 10^{-3}$ & 2.79 & $2.075 \cdot 10^{-3}$ & 2.79 \\
   & $64^3$ & $1.383 \cdot 10^{-3}$ & 2.93 & $5.435 \cdot 10^{-4}$ & 2.93 & $2.718 \cdot 10^{-4}$ & 2.93 \\
   & $128^3$ & $1.751 \cdot 10^{-4}$ & 2.98 & $6.878 \cdot 10^{-5}$ & 2.98 & $3.439 \cdot 10^{-5}$ & 2.98 \\
\midrule
  $B_2$ & $16^3$ & $1.514 \cdot 10^{-1}$ & --- & $5.910 \cdot 10^{-2}$ & --- & $2.912 \cdot 10^{-2}$ & --- \\
   & $32^3$ & $2.147 \cdot 10^{-2}$ & 2.82 & $8.418 \cdot 10^{-3}$ & 2.81 & $4.191 \cdot 10^{-3}$ & 2.80 \\
   & $64^3$ & $2.777 \cdot 10^{-3}$ & 2.95 & $1.090 \cdot 10^{-3}$ & 2.95 & $5.446 \cdot 10^{-4}$ & 2.94 \\
   & $128^3$ & $3.509 \cdot 10^{-4}$ & 2.98 & $1.378 \cdot 10^{-4}$ & 2.98 & $6.888 \cdot 10^{-5}$ & 2.98 \\
\midrule
  $B_3$ & $16^3$ & $1.889 \cdot 10^{-1}$ & --- & $7.373 \cdot 10^{-2}$ & --- & $3.634 \cdot 10^{-2}$ & --- \\
   & $32^3$ & $2.678 \cdot 10^{-2}$ & 2.82 & $1.050 \cdot 10^{-2}$ & 2.81 & $5.228 \cdot 10^{-3}$ & 2.80 \\
   & $64^3$ & $3.470 \cdot 10^{-3}$ & 2.95 & $1.362 \cdot 10^{-3}$ & 2.95 & $6.804 \cdot 10^{-4}$ & 2.94 \\
   & $128^3$ & $4.385 \cdot 10^{-4}$ & 2.98 & $1.722 \cdot 10^{-4}$ & 2.98 & $8.607 \cdot 10^{-5}$ & 2.98 \\
\midrule
  $E_1$ & $16^3$ & $2.184 \cdot 10^{-1}$ & --- & $8.632 \cdot 10^{-2}$ & --- & $4.314 \cdot 10^{-2}$ & --- \\
   & $32^3$ & $3.161 \cdot 10^{-2}$ & 2.79 & $1.245 \cdot 10^{-2}$ & 2.79 & $6.226 \cdot 10^{-3}$ & 2.79 \\
   & $64^3$ & $4.150 \cdot 10^{-3}$ & 2.93 & $1.631 \cdot 10^{-3}$ & 2.93 & $8.153 \cdot 10^{-4}$ & 2.93 \\
   & $128^3$ & $5.253 \cdot 10^{-4}$ & 2.98 & $2.063 \cdot 10^{-4}$ & 2.98 & $1.032 \cdot 10^{-4}$ & 2.98 \\
\midrule
  $E_2$ & $16^3$ & $1.129 \cdot 10^{-1}$ & --- & $4.409 \cdot 10^{-2}$ & --- & $2.174 \cdot 10^{-2}$ & --- \\
   & $32^3$ & $1.601 \cdot 10^{-2}$ & 2.82 & $6.278 \cdot 10^{-3}$ & 2.81 & $3.125 \cdot 10^{-3}$ & 2.80 \\
   & $64^3$ & $2.080 \cdot 10^{-3}$ & 2.94 & $8.167 \cdot 10^{-4}$ & 2.94 & $4.079 \cdot 10^{-4}$ & 2.94 \\
   & $128^3$ & $2.629 \cdot 10^{-4}$ & 2.98 & $1.032 \cdot 10^{-4}$ & 2.98 & $5.161 \cdot 10^{-5}$ & 2.98 \\
\midrule
  $E_3$ & $16^3$ & $7.153 \cdot 10^{-4}$ & --- & $2.791 \cdot 10^{-4}$ & --- & $1.370 \cdot 10^{-4}$ & --- \\
   & $32^3$ & $8.019 \cdot 10^{-5}$ & 3.16 & $3.151 \cdot 10^{-5}$ & 3.15 & $1.575 \cdot 10^{-5}$ & 3.12 \\
   & $64^3$ & $2.326 \cdot 10^{-6}$ & 5.11 & $9.140 \cdot 10^{-7}$ & 5.11 & $4.570 \cdot 10^{-7}$ & 5.11 \\
   & $128^3$ & $2.076 \cdot 10^{-7}$ & 3.49 & $8.154 \cdot 10^{-8}$ & 3.49 & $4.076 \cdot 10^{-8}$ & 3.49 \\
\bottomrule
\end{tabular}}
\caption{Errors of the face averages and experimental orders of convergence for the divergence-free method of
Section~\ref{sec:divfree} with initialization from the vector potentials on the oblique plane wave (\ref{eq:oblique_wave}) at $t=1$ and $\text{CFL} = 0.45$.}
\label{tab:eoc_oblique_divfree}
\end{table}

\begin{table}[htbp]
\centering
\small
\setlength{\tabcolsep}{4pt}
\adjustbox{max width=\linewidth}{%
\begin{tabular}{lrcrcrc}
\toprule
\textbf{$N^3$} & \textbf{classical} & \textbf{EOC} & \textbf{div-free, Simpson} & \textbf{EOC} & \textbf{div-free, potential} & \textbf{EOC} \\
\midrule
  $16^3$ & $4.173 \cdot 10^{-1}$ & --- & $1.095 \cdot 10^{-4}$ & --- & $9.804 \cdot 10^{-15}$ & --- \\
  $32^3$ & $1.124 \cdot 10^{-1}$ & 1.89 & $7.147 \cdot 10^{-6}$ & 3.94 & $2.376 \cdot 10^{-14}$ & --- \\
  $64^3$ & $2.451 \cdot 10^{-2}$ & 2.20 & $4.515 \cdot 10^{-7}$ & 3.98 & $1.308 \cdot 10^{-13}$ & --- \\
  $128^3$ & $6.256 \cdot 10^{-3}$ & 1.97 & $2.830 \cdot 10^{-8}$ & 4.00 & $4.375 \cdot 10^{-13}$ & --- \\
\bottomrule
\end{tabular}}
\caption{Divergence $\max\vert\nabla\cdot\bm B\vert$ of the reconstruction at $t=1$, with
$\text{CFL} = 0.45$ for all methods.}
\label{tab:eoc_div}
\end{table}

\begin{table}[htbp]
\centering
\small
\setlength{\tabcolsep}{4pt}
\adjustbox{max width=\linewidth}{%
\begin{tabular}{lrcrcrc}
\toprule
\textbf{$N^3$} & \textbf{classical} & \textbf{EOC} & \textbf{div-free, Simpson} & \textbf{EOC} & \textbf{div-free, potential} & \textbf{EOC} \\
\midrule
  $16^3$ & $4.069 \cdot 10^{-2}$ & --- & $9.908 \cdot 10^{-2}$ & --- & $9.905 \cdot 10^{-2}$ & --- \\
  $32^3$ & $5.617 \cdot 10^{-3}$ & 2.86 & $1.410 \cdot 10^{-2}$ & 2.81 & $1.410 \cdot 10^{-2}$ & 2.81 \\
  $64^3$ & $7.199 \cdot 10^{-4}$ & 2.96 & $1.829 \cdot 10^{-3}$ & 2.95 & $1.829 \cdot 10^{-3}$ & 2.95 \\
  $128^3$ & $9.065 \cdot 10^{-5}$ & 2.99 & $2.309 \cdot 10^{-4}$ & 2.99 & $2.309 \cdot 10^{-4}$ & 2.99 \\
\bottomrule
\end{tabular}}
\caption{Absolute value of the relative energy error (\ref{eq:energy_error}) at $t=1$, with
$\text{CFL} = 0.45$ for all methods.}
\label{tab:eoc_energy}
\end{table}

\subsection{Spherical discontinuity}
\begin{figure}[htbp]
	\centering
	\includegraphics[width=0.49\linewidth]{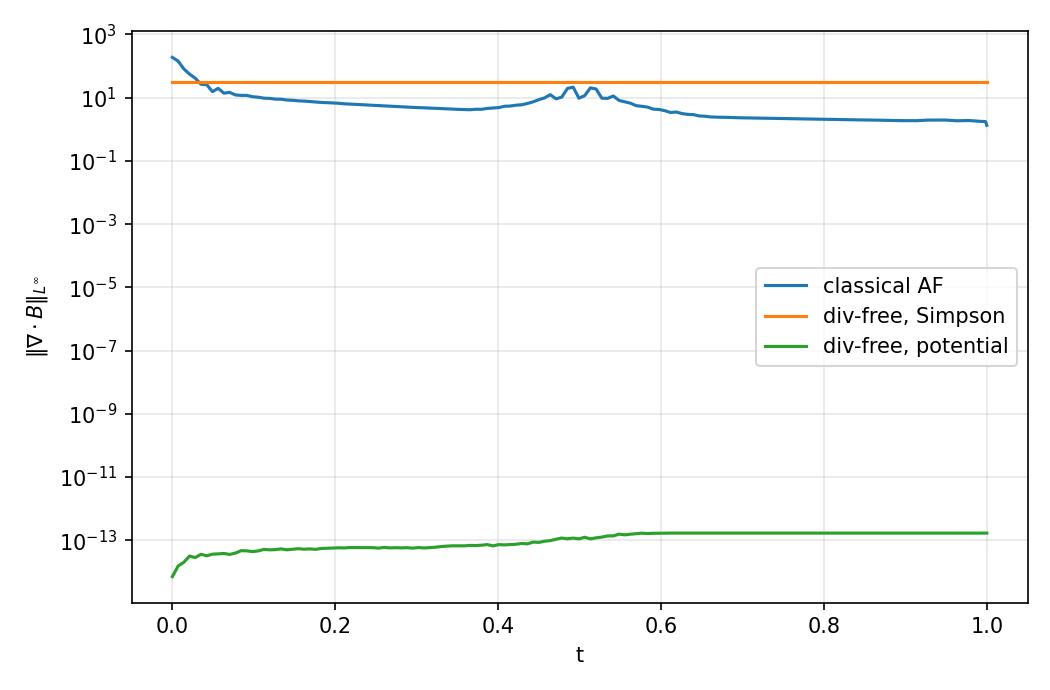}
	\includegraphics[width=0.49\linewidth]{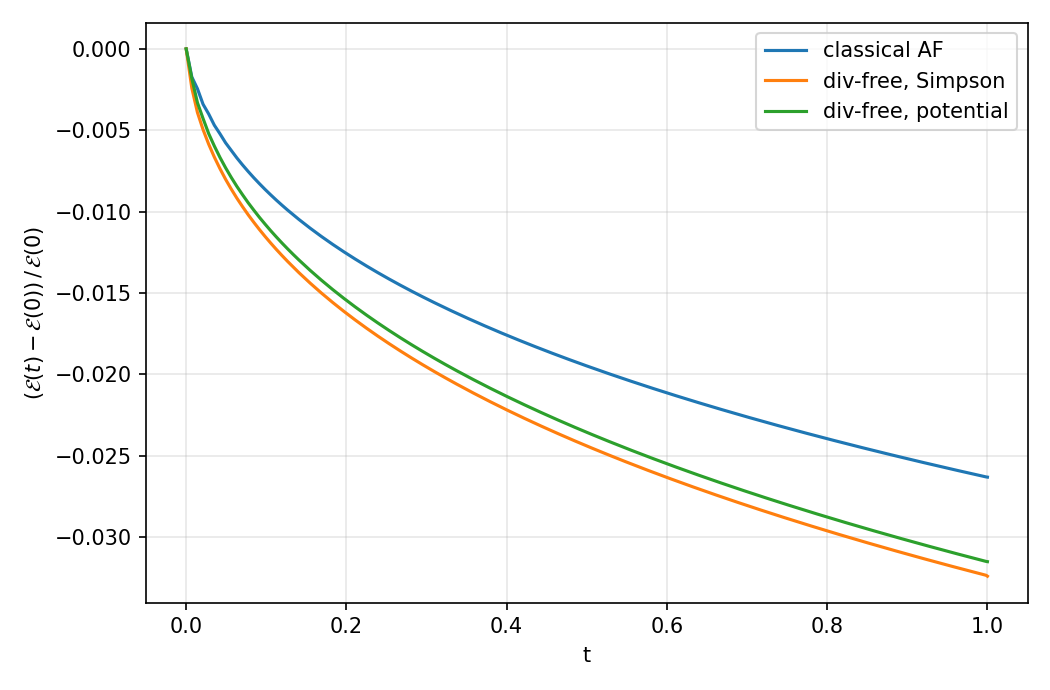}
	\caption{Maximum norm of $\nabla\cdot\bm B$ (left) and relative energy error (\ref{eq:energy_error}) (right)
	over $t \in [0,1]$ for example (\ref{eq:sphere}) on a grid with $128^3$ cells, for the classical Active Flux method and
	the divergence-free method with Simpson and potential initialization, with CFL $0.45$ for all methods.}
	\label{fig:discontinuity_divanden}
\end{figure}

\begin{figure}[htbp]
	\centering
	\includegraphics[width=\linewidth]{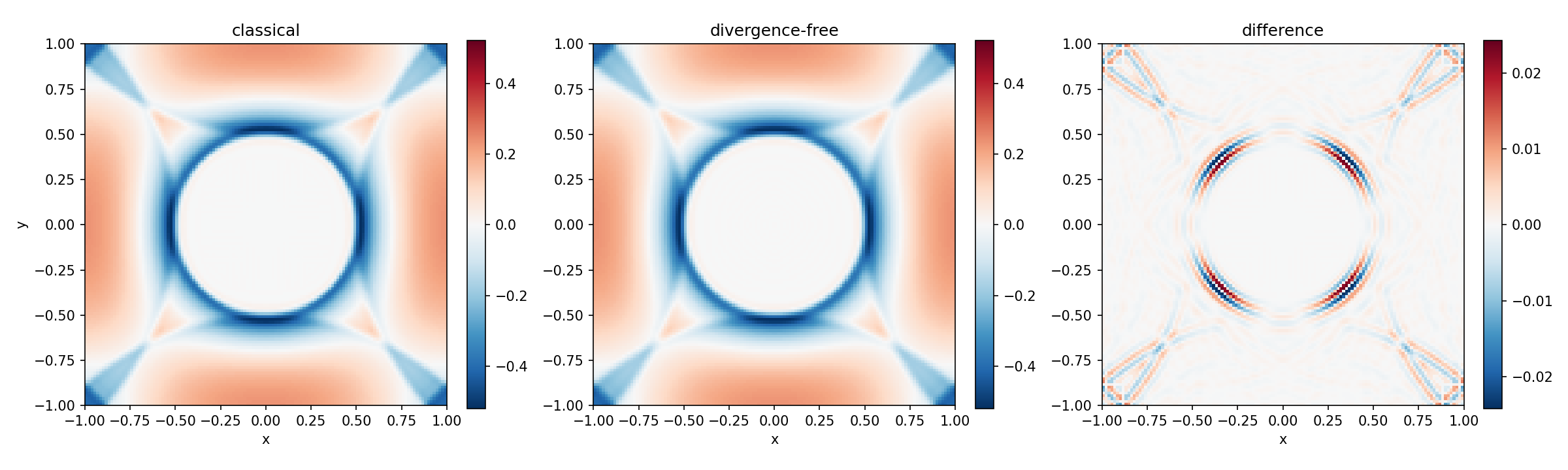}
	\caption{$E_3$ in the plane $z=0$ at $t = 1$ on a grid with $128^3$ cells and CFL $0.45$. Cell averages of the classical Active Flux method (left),
			 cell averages of the reconstruction of the divergence-free method with potential initialization (middle), difference of the two methods (right).}
	\label{fig:disc_plane}
\end{figure}

\begin{figure}[htbp]
	\centering
	\includegraphics[width=\linewidth]{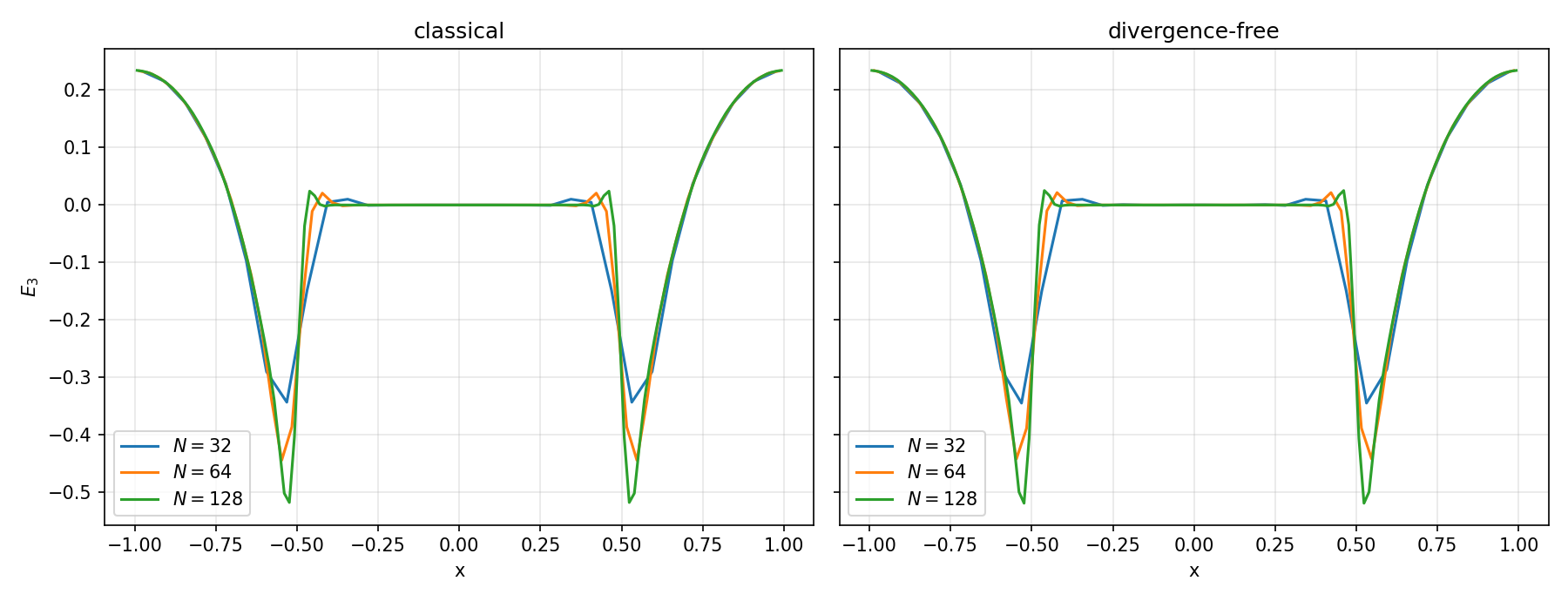}
	\caption{Cut of $E_3$ along the $x$-axis ($y=z=0$) at $t = 1$ for $N = 32, 64, 128$. Left: cell averages of the classical Active Flux method, right: cell averages of the
		reconstruction of the divergence-free method with potential initialization.}
	\label{fig:disc_slice}
\end{figure}

As a last example we want to consider a spherical discontinuity on the domain $[-1,1]^3$ with periodic boundary conditions.
The initial data is given by 
\begin{align}\label{eq:sphere}
	\bm B(\bm x, 0) = \begin{cases} \frac{1}{R_0}(-y, x, 0)^\top, & r < R_0, \\[4pt] \bm 0, & r \geq R_0, \end{cases}
	\qquad \bm E(\bm x, 0) = \bm 0,
\end{align}
with $r = \vert \bm x \vert$ and $R_0 = 0.5$.
We initialize the facial averages of $\bm B$ via Stokes' theorem from the potential
\begin{align}\label{eq:apot_sphere}
	\bm A(\bm x) = \begin{pmatrix}
		0 \\
		0 \\
		A_z(\bm x)
	\end{pmatrix}, \qquad
	A_z(\bm x) = -\frac{1}{2R_0}\begin{cases}
		x^2 + y^2, & r < R_0, \\
		R_0^2 - z^2, & r \geq R_0,\ \vert z \vert < R_0, \\
		0, & \vert z \vert \geq R_0.
	\end{cases}
\end{align}
In this example one discontinuity travels inwards and one outwards. The outgoing front reaches the boundary at time $t=0.5$, 
meaning that the solution at time $t=1$ in Figures~\ref{fig:disc_plane} and~\ref{fig:disc_slice} also contains the waves that have re-entered. 
In Figure~\ref{fig:discontinuity_divanden} we present the divergence and the relative energy error
for both methods. For the divergence-free method we again consider the two different initializations of the facial averages.
We see that even with the discontinuity the divergence-free method keeps the divergence at round-off if the face averages are initialized from the 
potential. Again, the Simpson initialization shows larger values, but these are preserved in time. 

In Figure~\ref{fig:disc_plane} we present $E_3$ in the plane $z = 0$ at $t = 1$ for $N^3 = 128^3$ cells for the classical 
and the divergence-free method. Besides that, we present the difference of the two planes. 
Furthermore, we present in Figure~\ref{fig:disc_slice} the respective cut at $y = 0$ for $N = 32, 64, 128$ for both methods. 
Both methods present oscillations at the discontinuities, but no method is superior to the other. 

\clearpage
\section{Conclusion}\label{sec:conclusion}
In this article we presented the first fully discrete three-dimensional Active Flux method 
for Maxwell's equations. We also presented an alternative approach towards an Active Flux 
method for Maxwell's equations that is divergence-free by design. Both methods are based on 
an evolution operator derived from Kirchhoff's formula. 
For the divergence-free method, the curl in the evolution operator is evaluated in the
distributional sense to account for the discontinuities of the reconstruction across cell faces.
In numerical examples, we compared accuracy, divergence preservation and energy error of 
both methods. Both methods are third-order accurate. The classical method runs stably up to
CFL $0.5$ and keeps the divergence error bounded. The divergence-free method preserves the 
divergence to round-off when the face averages are initialized from a vector potential and 
appears to be stable up to CFL $0.48$. Both methods dissipate energy, the classical method less than 
the divergence-free method; the dissipation of the divergence-free method is independent of 
the initialization of the face averages. For a spherical discontinuity, both methods show 
oscillations of comparable size.
In future work, we intend to incorporate source terms in order to couple these methods to 
Active Flux methods for the Vlasov equation.

\backmatter





\begin{appendices}

\section{Full description of the divergence-free reconstruction}\label{app:reconstruction}
With a cyclic permutation of the coefficients $a \to b, b \to c, c \to a, x \to y, y \to z, z \to x, \xi \to \eta, \eta \to \nu, \nu \to \xi$ in (\ref{eq:face_poly}) and (\ref{eq:divfree_ansatz}) 
we obtain the face reconstructions of $B_2$ and $B_3$ from the reconstruction of $B_1$. We have 
\begin{align*}
B_2(\xi, -1,\nu) &= \mathcal{B}^{y,-}_{0} P_0(\xi) P_0(\nu) + \mathcal{B}^{y,-}_{x} P_1(\xi) P_0(\nu) + \mathcal{B}^{y,-}_{z} P_0(\xi) P_1(\nu) \\
	&\phantom{=} + \mathcal{B}^{y,-}_{xx} P_2(\xi) P_0(\nu) + \mathcal{B}^{y,-}_{xz} P_1(\xi) P_1(\nu) + \mathcal{B}^{y,-}_{zz} P_0(\xi) P_2(\nu), \\
B_3(\xi, \eta,-1) &= \mathcal{B}^{z,-}_{0} P_0(\xi) P_0(\eta) + \mathcal{B}^{z,-}_{x} P_1(\xi) P_0(\eta) + \mathcal{B}^{z,-}_{y} P_0(\xi) P_1(\eta) \\
	&\phantom{=} + \mathcal{B}^{z,-}_{xx} P_2(\xi) P_0(\eta) + \mathcal{B}^{z,-}_{xy} P_1(\xi) P_1(\eta) + \mathcal{B}^{z,-}_{yy} P_0(\xi) P_2(\eta).
\end{align*}
And the respective cell reconstructions are given by
\begin{align*}
B_2^{rec}(\xi, \eta, \nu) = \ & b_0 + b_x \xi + b_y \eta + b_z \nu
	+ b_{xx} P_2(\xi) + b_{yy} P_2(\eta) + b_{zz} P_2(\nu) \\
	& + b_{xy} \xi\eta + b_{xz} \xi\nu + b_{yz} \eta\nu
	+ b_{xyz} \xi\eta\nu \\
	& + b_{yzz} \eta P_2(\nu) + b_{xxy} P_2(\xi) \eta
	+ b_{yyz} P_2(\eta) \nu \\
	& + b_{xyy} \xi P_2(\eta)
	+ b_{yyy} P_3(\eta), \\
B_3^{rec}(\xi, \eta, \nu) = \ & c_0 + c_x \xi + c_y \eta + c_z \nu
	+ c_{xx} P_2(\xi) + c_{yy} P_2(\eta) + c_{zz} P_2(\nu) \\
	& + c_{xy} \xi\eta + c_{xz} \xi\nu + c_{yz} \eta\nu
	+ c_{xyz} \xi\eta\nu \\
	& + c_{xxz} P_2(\xi) \nu + c_{yyz} P_2(\eta) \nu
	+ c_{xzz} \xi P_2(\nu) \\
	& + c_{yzz} \eta P_2(\nu)
	+ c_{zzz} P_3(\nu).
\end{align*}

\section{Coefficients of the divergence-free reconstruction}\label{app:coefficients}
Here we list the coefficients of the divergence-free reconstruction (\ref{eq:divfree_ansatz}) in terms of the facial degrees of freedom (\ref{eq:facial_FD}).
Those are derived by rescaling the reconstruction of \cite{BALSARA2009} onto the cube $[-1,1]^3$.
The coefficients of $B_1$ are given by:
\begin{align*}
a_{0}    &= \frac{1}{2}\bigl(\Bp{x}{0}+\Bm{x}{0}\bigr)
            +\frac{1}{6}\bigl(\Bp{y}{x}-\Bm{y}{x}\bigr)
            +\frac{1}{6}\bigl(\Bp{z}{x}-\Bm{z}{x}\bigr) \\
a_{x}    &= \frac{1}{2}\bigl(\Bp{x}{0}-\Bm{x}{0}\bigr)
            +\frac{1}{15}\bigl(\Bp{y}{xx}-\Bm{y}{xx}\bigr)
            +\frac{1}{15}\bigl(\Bp{z}{xx}-\Bm{z}{xx}\bigr) \\
a_{y}    &= \frac{1}{2}\bigl(\Bp{x}{y}+\Bm{x}{y}\bigr)
            +\frac{1}{12}\bigl(\Bp{z}{xy}-\Bm{z}{xy}\bigr) \\
a_{z}    &= \frac{1}{2}\bigl(\Bp{x}{z}+\Bm{x}{z}\bigr)
            +\frac{1}{12}\bigl(\Bp{y}{xz}-\Bm{y}{xz}\bigr) \\
a_{xx}   &= -\frac{1}{4}\bigl(\Bp{y}{x}-\Bm{y}{x}\bigr)
            -\frac{1}{4}\bigl(\Bp{z}{x}-\Bm{z}{x}\bigr) \\
a_{xy}   &= \frac{1}{2}\bigl(\Bp{x}{y}-\Bm{x}{y}\bigr) \\
a_{xz}   &= \frac{1}{2}\bigl(\Bp{x}{z}-\Bm{x}{z}\bigr) \\
a_{yy}   &= \frac{1}{2}\bigl(\Bp{x}{yy}+\Bm{x}{yy}\bigr) \\
a_{xyy}  &= \frac{1}{2}\bigl(\Bp{x}{yy}-\Bm{x}{yy}\bigr) \\
a_{zz}   &= \frac{1}{2}\bigl(\Bp{x}{zz}+\Bm{x}{zz}\bigr) \\
a_{xzz}  &= \frac{1}{2}\bigl(\Bp{x}{zz}-\Bm{x}{zz}\bigr) \\
a_{yz}   &= \frac{1}{2}\bigl(\Bp{x}{yz}+\Bm{x}{yz}\bigr) \\
a_{xyz}  &= \frac{1}{2}\bigl(\Bp{x}{yz}-\Bm{x}{yz}\bigr) \\
a_{xxx}  &= -\frac{1}{6}\bigl(\Bp{y}{xx}-\Bm{y}{xx}\bigr)
            -\frac{1}{6}\bigl(\Bp{z}{xx}-\Bm{z}{xx}\bigr) \\
a_{xxy}  &= -\frac{1}{8}\bigl(\Bp{z}{xy}-\Bm{z}{xy}\bigr) \\
a_{xxz}  &= -\frac{1}{8}\bigl(\Bp{y}{xz}-\Bm{y}{xz}\bigr).
\end{align*}
The coefficients of $B_2$ are given by:
\begin{align*}
b_{0}    &= \frac{1}{2}\bigl(\Bp{y}{0}+\Bm{y}{0}\bigr)
            +\frac{1}{6}\bigl(\Bp{x}{y}-\Bm{x}{y}\bigr)
            +\frac{1}{6}\bigl(\Bp{z}{y}-\Bm{z}{y}\bigr) \\
b_{y}    &= \frac{1}{2}\bigl(\Bp{y}{0}-\Bm{y}{0}\bigr)
            +\frac{1}{15}\bigl(\Bp{x}{yy}-\Bm{x}{yy}\bigr)
            +\frac{1}{15}\bigl(\Bp{z}{yy}-\Bm{z}{yy}\bigr) \\
b_{z}    &= \frac{1}{2}\bigl(\Bp{y}{z}+\Bm{y}{z}\bigr)
            +\frac{1}{12}\bigl(\Bp{x}{yz}-\Bm{x}{yz}\bigr) \\
b_{x}    &= \frac{1}{2}\bigl(\Bp{y}{x}+\Bm{y}{x}\bigr)
            +\frac{1}{12}\bigl(\Bp{z}{xy}-\Bm{z}{xy}\bigr) \\
b_{yy}   &= -\frac{1}{4}\bigl(\Bp{z}{y}-\Bm{z}{y}\bigr)
            -\frac{1}{4}\bigl(\Bp{x}{y}-\Bm{x}{y}\bigr) \\
b_{yz}   &= \frac{1}{2}\bigl(\Bp{y}{z}-\Bm{y}{z}\bigr) \\
b_{xy}   &= \frac{1}{2}\bigl(\Bp{y}{x}-\Bm{y}{x}\bigr) \\
b_{zz}   &= \frac{1}{2}\bigl(\Bp{y}{zz}+\Bm{y}{zz}\bigr) \\
b_{yzz}  &= \frac{1}{2}\bigl(\Bp{y}{zz}-\Bm{y}{zz}\bigr) \\
b_{xx}   &= \frac{1}{2}\bigl(\Bp{y}{xx}+\Bm{y}{xx}\bigr) \\
b_{xxy}  &= \frac{1}{2}\bigl(\Bp{y}{xx}-\Bm{y}{xx}\bigr) \\
b_{xz}   &= \frac{1}{2}\bigl(\Bp{y}{xz}+\Bm{y}{xz}\bigr) \\
b_{xyz}  &= \frac{1}{2}\bigl(\Bp{y}{xz}-\Bm{y}{xz}\bigr) \\
b_{yyy}  &= -\frac{1}{6}\bigl(\Bp{z}{yy}-\Bm{z}{yy}\bigr)
            -\frac{1}{6}\bigl(\Bp{x}{yy}-\Bm{x}{yy}\bigr) \\
b_{yyz}  &= -\frac{1}{8}\bigl(\Bp{x}{yz}-\Bm{x}{yz}\bigr) \\
b_{xyy}  &= -\frac{1}{8}\bigl(\Bp{z}{xy}-\Bm{z}{xy}\bigr).
\end{align*}
The coefficients of $B_3$ are given by:
\begin{align*}
c_{0}    &= \frac{1}{2}\bigl(\Bp{z}{0}+\Bm{z}{0}\bigr)
            +\frac{1}{6}\bigl(\Bp{x}{z}-\Bm{x}{z}\bigr)
            +\frac{1}{6}\bigl(\Bp{y}{z}-\Bm{y}{z}\bigr) \\
c_{z}    &= \frac{1}{2}\bigl(\Bp{z}{0}-\Bm{z}{0}\bigr)
            +\frac{1}{15}\bigl(\Bp{x}{zz}-\Bm{x}{zz}\bigr)
            +\frac{1}{15}\bigl(\Bp{y}{zz}-\Bm{y}{zz}\bigr) \\
c_{x}    &= \frac{1}{2}\bigl(\Bp{z}{x}+\Bm{z}{x}\bigr)
            +\frac{1}{12}\bigl(\Bp{y}{xz}-\Bm{y}{xz}\bigr) \\
c_{y}    &= \frac{1}{2}\bigl(\Bp{z}{y}+\Bm{z}{y}\bigr)
            +\frac{1}{12}\bigl(\Bp{x}{yz}-\Bm{x}{yz}\bigr) \\
c_{zz}   &= -\frac{1}{4}\bigl(\Bp{x}{z}-\Bm{x}{z}\bigr)
            -\frac{1}{4}\bigl(\Bp{y}{z}-\Bm{y}{z}\bigr) \\
c_{xz}   &= \frac{1}{2}\bigl(\Bp{z}{x}-\Bm{z}{x}\bigr) \\
c_{yz}   &= \frac{1}{2}\bigl(\Bp{z}{y}-\Bm{z}{y}\bigr) \\
c_{xx}   &= \frac{1}{2}\bigl(\Bp{z}{xx}+\Bm{z}{xx}\bigr) \\
c_{xxz}  &= \frac{1}{2}\bigl(\Bp{z}{xx}-\Bm{z}{xx}\bigr) \\
c_{yy}   &= \frac{1}{2}\bigl(\Bp{z}{yy}+\Bm{z}{yy}\bigr) \\
c_{yyz}  &= \frac{1}{2}\bigl(\Bp{z}{yy}-\Bm{z}{yy}\bigr) \\
c_{xy}   &= \frac{1}{2}\bigl(\Bp{z}{xy}+\Bm{z}{xy}\bigr) \\
c_{xyz}  &= \frac{1}{2}\bigl(\Bp{z}{xy}-\Bm{z}{xy}\bigr) \\
c_{zzz}  &= -\frac{1}{6}\bigl(\Bp{x}{zz}-\Bm{x}{zz}\bigr)
            -\frac{1}{6}\bigl(\Bp{y}{zz}-\Bm{y}{zz}\bigr) \\
c_{xzz}  &= -\frac{1}{8}\bigl(\Bp{y}{xz}-\Bm{y}{xz}\bigr) \\
c_{yzz}  &= -\frac{1}{8}\bigl(\Bp{x}{yz}-\Bm{x}{yz}\bigr).
\end{align*}

\end{appendices}

\bibliography{bibliography}

\end{document}